\documentclass[11pt,a4paper]{article}
\usepackage[T1]{fontenc}
\usepackage[utf8]{inputenc}
\usepackage{lmodern}
\usepackage[top=28mm,bottom=28mm,left=29mm,right=29mm,headheight=14pt]{geometry}
\usepackage{amsmath,amssymb,amsthm,mathtools}
\usepackage[final]{microtype}
\usepackage{fancyhdr}
\usepackage{xcolor}
\usepackage{needspace,etoolbox}
\usepackage[colorlinks=true,linkcolor=blue!45!black,
  citecolor=green!35!black,urlcolor=blue!55!black]{hyperref}
\numberwithin{equation}{section}
\hypersetup{pdftitle={Quantum wedges and large deviations for SLE with two force points}}
\makeatletter
\renewcommand{\@maketitle}{%
  \newpage\null\vskip0.5em
  \begin{center}
    {\LARGE\@title\par}
    \ifx\@author\@empty\else
      \vskip1em{\large\begin{tabular}[t]{c}\@author\end{tabular}\par}
    \fi
    \ifx\@date\@empty\else\vskip0.75em{\@date\par}\fi
  \end{center}
  \vskip0.75em}
\makeatother
\newtheorem{theorem}{Theorem}[section]
\newtheorem{proposition}[theorem]{Proposition}
\newtheorem{lemma}[theorem]{Lemma}
\newtheorem{corollary}[theorem]{Corollary}
\theoremstyle{definition}
\newtheorem{definition}[theorem]{Definition}
\theoremstyle{remark}
\newtheorem{remark}[theorem]{Remark}
\BeforeBeginEnvironment{theorem}{\Needspace{7\baselineskip}}
\BeforeBeginEnvironment{proposition}{\Needspace{7\baselineskip}}
\BeforeBeginEnvironment{lemma}{\Needspace{7\baselineskip}}
\BeforeBeginEnvironment{corollary}{\Needspace{7\baselineskip}}
\newcommand{\Sstrip}{\mathcal S}
\newcommand{\R}{\mathbb R}
\newcommand{\D}{\mathbb D}
\newcommand{\E}{\mathbb E}
\newcommand{\Pp}{\mathbb P}
\newcommand{\dd}{\,d}
\DeclareMathOperator{\dist}{dist}
\DeclareMathOperator{\Rea}{Re}
\DeclareMathOperator{\Ima}{Im}
\DeclareMathOperator{\supp}{supp}
\title{Quantum wedges and large deviations\\
for SLE with two force points}
\author{
  \begin{tabular}{@{}c@{\qquad}c@{}}
    Yuchen Fan\thanks{These authors contributed equally
      and are joint first authors.}
    & Zhenfeng Tu\footnotemark[1] \\[0.4em]
    {\small Courant Institute}
    & {\small Courant Institute} \\
    {\small New York University}
    & {\small New York University} \\
    {\small\href{mailto:yf2085@nyu.edu}{yf2085@nyu.edu}}
    & {\small\href{mailto:zt2255@nyu.edu}{zt2255@nyu.edu}}
  \end{tabular}%
}
\date{}
\begin{document}
\raggedbottom
\maketitle

\begin{abstract}
We prove a good large deviation principle for chordal
$\mathrm{SLE}_\kappa(W_1-2;W_2-2)$ with force points at $0^-$ and
$0^+$, for every fixed $W_1,W_2>0$, as $\kappa\downarrow0$.
The topology is Hausdorff convergence of the entire compact trace.
The rate has two equivalent descriptions: the sum of the component
wedge energies obtained by cutting a deterministic ambient wedge,
and the Loewner driving energy relative to the force-point drift.
\end{abstract}

\section{Introduction}\label{sec:introduction}

Schramm--Loewner evolution (SLE), introduced by Schramm \cite{Schramm},
is a family of conformally invariant random curves arising as scaling
limits of interfaces in two-dimensional statistical mechanics. Chordal
$\mathrm{SLE}_\kappa$ in the upper half-plane is generated by
\[
 \partial_tg_t(z)=\frac{2}{g_t(z)-\sqrt\kappa B_t},\qquad g_0(z)=z,
\]
where $B$ is standard Brownian motion. The map $g_t$ is the conformal
map from the complement of the Loewner hull $K_t$ onto $\mathbb H$,
normalized by $g_t(z)=z+2t/z+O(z^{-2})$ at infinity. The hull contains
the trace up to time $t$ and any regions it has disconnected from
infinity; the trace itself is recovered from the inverse maps.

SLE with force points incorporates additional marked boundary points
into the driving equation \cite{SW,IG1}. We use the standard
force-parameter convention:
\begin{definition}[SLE with two boundary force points]
\label{def:SLE}
For $0<\kappa<4$ and $\rho_1,\rho_2>-2$, chordal
$\mathrm{SLE}_\kappa(\rho_1;\rho_2)$ from $0$ to $\infty$ in
$\mathbb H$, with force points at $0^-,0^+$, is the trace generated by
\begin{equation}\label{eq:SLE}
 \begin{aligned}
 dU_t&=\sqrt\kappa\,dB_t
       +\frac{\rho_1}{U_t-V_t^-}\,dt
       +\frac{\rho_2}{U_t-V_t^+}\,dt,\\
 dV_t^\pm&=\frac{2}{V_t^\pm-U_t}\,dt,\\
 \partial_tg_t(z)&=\frac{2}{g_t(z)-U_t},\qquad g_0(z)=z.
 \end{aligned}
\end{equation}
Here $V_t^\pm$ are the Loewner images of the force points. The initial
condition $U_0=V_0^-=V_0^+=0$ is interpreted with the ordered force
points $0^-,0^+$, using the standard continuation of the singular
system \cite[Section 2.2]{IG1}. The trace is
\[
 \eta(t)=\lim_{y\downarrow0}g_t^{-1}(U_t+iy),\qquad t\ge0.
\]
\end{definition}

A large deviation principle (LDP) describes rare events on an
exponential scale \cite{DZ}. For a family of random objects indexed by
$\kappa\downarrow0$, an LDP at speed $\kappa^{-1}$ assigns a
nonnegative cost $I(x)$ to each configuration $x$:
informally, observing a configuration near $x$ has probability
on the scale $\exp(-I(x)/\kappa)$. The precise formulation is
as follows.

\begin{definition}[Large deviation principle]
Let $(X_\kappa)_{\kappa>0}$ be random variables taking values in a metric space $E$, and let $I:E\to[0,\infty]$ be lower semicontinuous. We say that $(X_\kappa)$ satisfies an LDP at speed $\kappa^{-1}$ with rate function $I$ if
\[
\begin{aligned}
\liminf_{\kappa\downarrow0}
\kappa\log\mathbb P(X_\kappa\in G)
&\geq -\inf_{x\in G}I(x)
&&\text{for every open }G\subset E,\\
\limsup_{\kappa\downarrow0}
\kappa\log\mathbb P(X_\kappa\in F)
&\leq -\inf_{x\in F}I(x)
&&\text{for every closed }F\subset E.
\end{aligned}
\]
Here $\inf\varnothing=+\infty$ and $\log0=-\infty$. The rate function is called \emph{good} if its sublevel sets $\{x\in E:I(x)\leq a\}$ are compact for every finite $a\geq0$.
\end{definition}

Large deviations for chordal and multichordal SLE were established
by Wang \cite{WangLDP} and Peltola--Wang \cite{PW}, using the driving
functions; see Wang's survey \cite{WangSurvey} for an overview of
these results and their connections with conformal geometry.
Krusell \cite{Krusell} proved large deviations with one
boundary force point in the chordal and radial settings. Guskov \cite{Guskov2023} established a chordal LDP for capacity-parameterized curves in the topology of locally uniform convergence in time. Abuzaid--Peltola \cite{AbuzaidPeltola2026} obtained LDPs for full chordal, radial, and multichordal curves, including their endpoints, and for unparameterized curves. We study
\[
\mathrm{SLE}_\kappa(W_1-2;W_2-2),
\qquad W_1,W_2>0,
\]
as $\kappa\downarrow0$, with the weights fixed. The coincident force
points give a singular initial driving equation. We approach the
curve LDP through quantum welding, motivated by the relation between
Loewner energy, Dirichlet energy, and conformal welding in \cite{VW}.

The starting point is the quantum wedge welding identity
\cite[Theorems~1.2 and~1.4]{DMS}, also stated in \cite[Theorem~2.26]{AHS}. This welding construction builds on Sheffield's quantum zipper
\cite{Sheffield2016Zipper}.Write
\[
\gamma=\sqrt{\kappa},\qquad
\rho_i=W_i-2,\qquad W=W_1+W_2.
\]
We use these weight conventions throughout.
For $\kappa<\min\{4,2W_1,2W_2\}$, the interface is almost surely
simple and has no real boundary contacts other than its endpoints
\cite{IG1}.
For all sufficiently small $\kappa$, welding two independent
quantum wedges of weights $W_1,W_2$ along their boundary length
measures produces a weight-$W$ wedge decorated with an independent
$\mathrm{SLE}_\kappa(\rho_1;\rho_2)$ interface. Conversely, cutting
a weight-$W$ wedge along such an independent interface produces
the two independent component wedges.

This identity suggests a corresponding balance of large deviation
costs:
\[
\text{component costs}
=
\text{ambient cost}+\text{interface cost}.
\]

Wang \cite[Section 3.1]{WangSurvey} explains the cutting--welding
energy identity through a heuristic large deviation argument for
the quantum zipper: independence turns the costs of the ambient
field and interface into the sum of the component field costs.
Our joint LDP and finite-rate welding stability make this principle
rigorous for quantum wedges with the fixed weights considered here.
This is distinct from the question in
\cite[Remark 2.13]{WangSurvey} concerning an LDP for SLE loop
measures, whose infinite total mass requires an appropriate
normalization. Our chordal probability laws do not address that
loop-measure question.

We now specify the curve topology and the finite-rate stability
properties needed for this argument.
\subsection{Curves and topology}

All quantum fields and their conformal embeddings are defined on
\[
 \mathcal S=\mathbb R\times(0,\pi).
\]
We use the disk only to define a bounded metric for curves and to
obtain compactness of the conformal maps. The fixed coordinate map is
\begin{equation}\label{eq:chart}
 \chi(z)=\frac{e^z-i}{e^z+i}:\mathcal S\longrightarrow\mathbb D,
 \qquad \mathbb D=\{z\in\mathbb C:|z|<1\}.
\end{equation}
Adjoin the two ends to the closed strip and write
\[
 \overline{\mathcal S}^{\,*}
 =\bigl(\mathbb R\times[0,\pi]\bigr)\cup\{-\infty,+\infty\}.
\]
The map $\chi$ extends to a homeomorphism
$\overline\chi:\overline{\mathcal S}^{\,*}\to\overline{\mathbb D}$,
with $\overline\chi(-\infty)=-1$ and
$\overline\chi(+\infty)=1$. The two boundary lines map to the
lower and upper semicircles, respectively. We equip this compactified
strip with the metric
\[
 d_\chi(z,w)=|\overline\chi(z)-\overline\chi(w)|.
\]
Thus the geometry remains in strip coordinates; the disk supplies
its metric, rather than a different embedding of the fields.

An interface is a simple curve in the strip joining its two ends.
More precisely, it is the image $\eta=c([0,1])$ of a continuous
injective map
\[
 c:[0,1]\longrightarrow\overline{\mathcal S}^{\,*},\qquad
 c(0)=+\infty,\quad c(1)=-\infty,\quad c((0,1))\subset\mathcal S.
\]
We retain the curve including its two ends, but not its
parametrization. Let $\mathcal X$ denote the set of these curves.
For $\eta\in\mathcal X$, write
\[
 \widehat\eta=\overline\chi(\eta)
 =\overline{\chi(\eta\cap\mathcal S)}\subset\overline{\mathbb D}.
\]
The topology on $\mathcal X$ is defined by
\begin{equation}\label{eq:stripHausdorff}
 d_{\mathrm H}^{\chi}(\eta,\widetilde\eta)
 :=d_{\mathrm H}^{\overline{\mathbb D}}
       (\widehat\eta,\widehat{\widetilde\eta}),
\end{equation}
where $d_{\mathrm H}^{\overline{\mathbb D}}$ is Euclidean
Hausdorff distance in the closed disk. In particular, we study
curves in $\mathcal S$ with the topology obtained by mapping their
entire traces to the disk. This is not Euclidean Hausdorff distance
on the unbounded strip.

The distance $d_{\mathrm H}^{\chi}$ is not conformally invariant,
but its induced topology does not depend on the choice of the
conformal chart $\chi:\mathcal S\to\mathbb D$. Indeed, any other
such chart is $m\circ\chi$ for a disk automorphism $m$, and both
$m$ and $m^{-1}$ extend uniformly continuously to the closed disk.
Thus a change of chart preserves Hausdorff convergence, although
it generally changes numerical distances.

This topology should be distinguished from the Carath\'eodory
topology, which describes kernel convergence of complementary
domains through convergence of suitably normalized conformal maps
on interior compact sets. For general hulls, Hausdorff and
Carath\'eodory convergence do not imply one another; see
\cite[Sections~2.1--2.2]{PW}. In particular, interior map
convergence alone does not identify the Hausdorff limit of an
entire curve. Our proof therefore includes a separate identification
of that limit, including its possible contacts with the outer
boundary.

For a curve $\eta\in\mathcal X$, let $D_1,D_2$ be the components
of $\mathcal S\setminus\eta$ adjacent to $\theta=0,\pi$,
respectively. The geometric maps $G_i:\mathcal S\to D_i$ are the
conformal bijections preserving both ends and normalized by
\begin{equation}\label{eq:geomark}
 G_1(0)=0,\qquad G_2(i\pi)=i\pi.
\end{equation}
They are determined by the curve. These maps are used to define
the cut fields in the rate function below.

\subsection{The curve rate}

For a weight $w>0$ and normalized boundary $j$, define
\begin{equation}\label{eq:energy}
I_w^{(j)}(\phi)
=
\frac1{4\pi}\int_{\mathcal S}
|\nabla(\phi+w\operatorname{Re}z)|^2\,d^2z,
\qquad
\int_0^\infty e^{\phi(t,j)/2}\,dt=1.
\end{equation}
The value is infinite unless the normalization holds and the field
is locally $H^1$ with finite displayed energy. The unique zero is
\begin{equation}\label{eq:zero}
\phi_w^*(t,\theta)=-wt+2\log(w/2).
\end{equation}

For $\eta\in\mathcal X$, define the cut fields of the zero-energy
ambient profile by
\begin{equation}\label{eq:cutzero}
 \phi_i^\eta=\phi_W^*\circ G_i+2\log|G_i'|.
\end{equation}
Our main result is the following.
\begin{theorem}\label{thm:main}
Fix $W_1,W_2>0$. The strip interfaces of
$\mathrm{SLE}_\kappa(W_1-2;W_2-2)$ satisfy a good LDP at speed
$\kappa^{-1}$ on $(\mathcal X,d_{\mathrm H}^\chi)$, with rate
\begin{equation}\label{eq:mainrate}
J(\eta)
=
I_{W_1}^{(0)}(\phi_1^\eta)
+
I_{W_2}^{(\pi)}(\phi_2^\eta).
\end{equation}
The rate is infinite unless both cut fields have finite energy.
The same rate gives a good LDP for the pair of geometrically
normalized complementary maps, with local uniform convergence
in the metric $d_\chi$.
\end{theorem}

We also identify this rate with the force-point driving energy.
After mapping the curve to the upper half-plane, the identity is
\[
J(\eta)
=
\frac12\int_0^\infty
\left|
\dot U_t-\frac{W_1-2}{U_t-V_t^-}
         -\frac{W_2-2}{U_t-V_t^+}
\right|^2\,dt.
\]
This expression has value $+\infty$ unless the driver is locally
absolutely continuous on $(0,\infty)$. The integral is understood
as the increasing limit over compact positive time intervals.
The equality includes equality of the finite domains.

\subsection{Proof strategy}

The field LDP does not directly give the curve LDP, since welding is
not known to be continuous in the distributional field topology used
here. Instead, we identify limits of actual welding configurations
when the limiting fields have finite rate. We organize the proof using
three groups of data: $A$ records the
ambient field, its boundary measures, and cumulative
boundary lengths; $B$ records the corresponding data for the two
independent component wedges; and $M$ records the cutting maps
and interface. Their precise definitions and
state spaces are given in Section~\ref{sec:field}. The argument
has five steps.

\begin{enumerate}
 \item \textbf{Enhanced marginal LDPs.}
 We first compute the field rate and then adjoin the boundary
 measures and cumulative boundary lengths. The
 difficulty is that GMC is not continuous in the distributional
 field topology. We apply contraction to continuous regularized
 approximations and show that the approximation errors are
 negligible at speed $\kappa^{-1}$, following the exponential
 approximation approach in
 Fan--Tu \cite{FT}.
 Exponential boundary tail estimates allow us to pass from local
 measures to cumulative lengths and implement the chosen embeddings.
 This gives good marginal rates $I(A)$ and $K(B)$; independence
 makes $K$ the sum of the two component rates.

 \item \textbf{Subsequential joint LDPs.}
 We adjoin the maps and interface of the actual quantum welding
 configurations. The enhanced field marginals are exponentially
 tight, while the enlarged map and trace spaces are compact.
 Hence every sequence $\kappa_n\downarrow0$ has a subsequence
 with a good joint LDP for $(A,B,M)$. Write $L$ for its initially
 unknown rate. Marginal contraction gives
 \[
 I(A)=\inf_{B,M}L(A,B,M),\qquad
 K(B)=\inf_{A,M}L(A,B,M).
 \]
 These identities determine the minimum cost above each marginal
 datum, but do not yet determine the full joint rate.

\item \textbf{Stability of conformal welding at finite rate}

The subsequence construction does not yet identify the joint
rate function: different subsequences could, a priori, yield
different rates $L$. To determine $L$, we will show that, for each
$B$ with $K(B)<\infty$, there is at most one pair $(A,M)$ at which
the joint rate can be finite. Goodness then makes the marginal
infimum attained, and
\[
K(B)=\inf_{A,M}L(A,B,M)
\]
forces the cost of that configuration to equal $K(B)$.

If $L(z)<\infty$, the LDP lower bound allows us to select actual
configurations $z_n$ converging to $z=(A,B,M)$, and marginal
contraction gives $I(A)+K(B)<\infty$. We prove that every such limit
satisfies the classical cutting and boundary length identities.
The limiting maps fill the ambient strip, and their separating curve
equals the recorded Hausdorff limit. Classical welding uniqueness
then determines the full pair $(A,M)$ from $B$.

 \item \textbf{Projection and identification of the curve rate.}
 Projection gives a good LDP for the interface. The exact quantum
 welding identity makes the ambient wedge independent of the
 SLE curve, so their joint rate is $I(A)+J(\eta)$. Evaluating
 at the zero-energy ambient profile identifies $J$ with the sum
 of the two cut-field energies in \eqref{eq:mainrate}. Conversely,
 a curve with finite-energy cut fields supplies classical
 component data; uniqueness identifies its configuration with
 the one furnished by the joint LDP. This proves the stated
 finite domain as well as the rate formula.

 \item \textbf{Identification with the driving energy.}
 A deterministic Dirichlet projection argument identifies
 this conformal rate with the force-point driving energy. In
 particular, finite conformal energy implies local absolute
 continuity of the driver on positive time intervals; this
 regularity is derived rather than assumed.
\end{enumerate}

\subsubsection{Strategy for SLE Loops}
The conformal welding identity of Ang--Holden--Sun
\cite[Theorem 1.1]{AHSLoop} suggests a possible extension to the question of finding the LDP 
SLE loop-measure  in \cite[Remark 2.13]{WangSurvey}.
Writing $\mathrm{QS}_\gamma$ for their unmarked quantum sphere
measure and $\mathrm{QD}_\gamma(\ell)$ for their unmarked quantum
disk measure with boundary length $\ell$, the identity is
\[
 \mathrm{QS}_\gamma\otimes\mathrm{SLE}^{\mathrm{loop}}_\kappa
 =C_\gamma\int_0^\infty
 \ell\,\operatorname{Weld}\bigl(
 \mathrm{QD}_\gamma(\ell),\mathrm{QD}_\gamma(\ell)\bigr)\,d\ell,
 \qquad \kappa=\gamma^2.
\]
Here $C_\gamma>0$, and $\operatorname{Weld}$ includes a uniformly
chosen boundary identification. The disks are the unmarked versions
of weight-$2$ disks; the unmarked sphere is obtained from a
weight-$(4-\gamma^2)$ sphere by reweighting by inverse area squared
and removing its two area-sampled marks
\cite[Definitions 2.3 and 2.6]{AHSLoop}.

Following \cite[Section 3.1]{WangSurvey} and
\cite[Section 1.3]{VW}, the expected rate balance is schematically
\[
 \text{joint cost of the two disks}
 =\text{sphere cost}+\text{loop cost},
\]
with compatible normalizations. A possible approach is to:
\begin{enumerate}
 \item establish joint disk LDPs retaining both area and boundary
 length, building on the unit-area result of \cite{FT};
 \item establish a unit-area sphere LDP and impose the corresponding
 constraints on the disks: equal boundary lengths and areas summing
 to one;
 \item prove finite-rate welding stability and uniqueness for these
 data, then identify the loop rate by projection as in this paper.
\end{enumerate}
This remains a proposed strategy. The loop measure is infinite,
and fixing the sphere area alone does not normalize it. One must
specify a finite restriction or another suitable normalization,
control the length integral, and track the welding constant at
the large deviation scale.

\section*{Acknowledgements}

The authors thanks Prof. Nina Holden for useful discussions. The authors used AI for motivation and for polishing the presentation of the paper. The convergence to classical conformal welding part is completed with substantial help from AI. Even 3 months ago, the current result seems hanging high up on the sky and we only have some broad plans on what to prove. With the help of AI the process is accelerated much more. Man! What can we say!

\section{Quantum wedges, welding data, and field large deviations}\label{sec:field}

\subsection{Definition of quantum wedges}

Fix $\gamma\in(0,2)$ and write
\[
Q=\frac{2}{\gamma}+\frac{\gamma}{2},
\qquad \kappa=\gamma^2.
\]

\begin{definition}[Quantum surface]\label{def:quantumsurface}
A $\gamma$-quantum surface is an equivalence class of pairs
$(D,h)$, where $D$ is a planar domain and $h$ is a distribution
on $D$. Two representatives $(D,h)$ and
$(\widetilde D,\widetilde h)$ are equivalent if there is a
conformal bijection $F:\widetilde D\to D$ such that
\[
\widetilde h=h\circ F+Q\log|F'|.
\]
For surfaces with ordered marked points, we additionally require
$F$ to preserve those points. Boundary marked points are
interpreted as prime ends when necessary.

An \emph{embedding} is a choice of representative of the
equivalence class.
\end{definition}

This is the quantum surface convention of
\cite[Definition 1.1]{DMS}. For the scaled field $\phi=\gamma h$,
the coordinate rule becomes
\[
\widetilde\phi
=\phi\circ F+(2+\kappa/2)\log|F'|.
\]

We next describe the Gaussian field used in the strip;
for background on the Neumann GFF, see
\cite[Chapter~6]{BerestyckiPowell2025}. Equip smooth real-valued
functions of finite Dirichlet energy with the inner product
\[
\langle f,g\rangle_\nabla
=\frac{1}{2\pi}\int_{\mathcal S}
\nabla f\cdot\nabla g\,d^2z.
\]
Let $H(\mathcal S)$ be the Hilbert completion after fixing
additive constants by
\[
\int_0^\pi f(i\theta)\,d\theta=0.
\]
Its orthogonal decomposition is
\[
H(\mathcal S)
=H^{\mathrm{const}}(\mathcal S)\oplus H^\perp(\mathcal S),
\]
where $H^{\mathrm{const}}(\mathcal S)$ consists of functions
depending only on the horizontal coordinate, and
\[
H^\perp(\mathcal S)
=
\left\{
f\in H(\mathcal S):
\int_0^\pi f(t+i\theta)\,d\theta=0
\text{ for almost every }t
\right\}.
\]
Orthogonality follows by integrating the horizontal derivative
of a function in $H^\perp(\mathcal S)$ over each vertical segment.

Let $(e_n)_{n\ge1}$ be an orthonormal basis of
$H^\perp(\mathcal S)$ and let $(a_n)_{n\ge1}$ be independent
standard Gaussian variables. The \emph{lateral Neumann GFF} is
\[
X^\perp=\sum_{n\ge1}a_ne_n,
\]
with the series interpreted as a random distribution. Its law
does not depend on the choice of basis. Its vertical average
vanishes, and its boundary covariance is
$-2\log|t-u|+O(1)$ on bounded intervals.

\begin{definition}[Weight-$w$ quantum wedge]
\label{def:quantumwedge}
Let $w>\gamma^2/2$ and put
\[
a_\gamma=\frac{w}{\gamma}-\frac{\gamma}{2}>0.
\]
Let $R_\gamma$ be a continuous process with $R_\gamma(0)=0$
whose two halves are independent and satisfy:
\begin{enumerate}
\item For $t\ge0$,
\[
R_\gamma(t)=B_{2t}-a_\gamma t,
\]
where $B$ is standard Brownian motion.
\item The process $(R_\gamma(-t))_{t\ge0}$ is a variance-two
Brownian motion with drift $a_\gamma$, conditioned to stay
positive and started from its entrance law at zero.
\end{enumerate}
Independently sample $X^\perp$ and set
\[
h_\gamma(t+i\theta)=R_\gamma(t)+X^\perp(t+i\theta).
\]
The doubly marked quantum surface represented by
\[
(\mathcal S,h_\gamma;+\infty,-\infty)
\]
is a \emph{weight-$w$ quantum wedge}. This representative is
its standard strip embedding.
\end{definition}

This is the thick-wedge construction in
\cite[Section 2.3]{AHS}. For every fixed $w>0$, the condition
$w>\gamma^2/2$ holds for all sufficiently small $\gamma$, so this
construction covers the semiclassical regime considered here.

Write $\phi_\kappa=\gamma h_\gamma$. The boundary length
measures are defined by the usual subcritical GMC construction:
\[
\nu_{\kappa,j}(dt)
=
e^{\gamma R_\gamma(t)/2}\nu^\perp_{\kappa,j}(dt),
\qquad j\in\{0,\pi\},
\]
where $\nu^\perp_{\kappa,j}$ is the boundary GMC of $X^\perp$.
Their coordinate rule is
\[
\nu_{\widetilde h}(A)=\nu_h(F(A)),
\qquad
\widetilde h=h\circ F+Q\log|F'|.
\]
We use this geometric normalization throughout. On bounded boundary
intervals it differs from full Wick subtraction for the lateral
field by a smooth deterministic multiplier tending to one as
$\kappa\downarrow0$.

\begin{lemma}[Radial Gaussian representation]\label{lem:radial}
Let $\mathbf B$ be a three-dimensional Brownian motion whose
coordinates have variance $2t$, and let $e_1=(1,0,0)$. Then
\[
 \bigl(\gamma R_\gamma(-t)\bigr)_{t\ge0}
 \stackrel{\mathrm{law}}=
 \left(\left|\sqrt\kappa\mathbf B_t
             +(w-\kappa/2)t e_1\right|\right)_{t\ge0}.
\]
Independently, the positive-time scaled radial part is
$\sqrt\kappa B_{2t}-(w-\kappa/2)t$.
\end{lemma}
\begin{proof}
For $a>0$, a variance-two Brownian motion with drift $a$, killed
at zero, has generator $f''+af'$. Conditioning it to stay positive
is the Doob transform by $h_a(r)=1-e^{-ar}$, with generator
\[
 f''(r)+\left(a+2\frac{h_a'(r)}{h_a(r)}\right)f'(r)
 =f''(r)+a\coth(ar/2)f'(r).
\]
Its entrance law at zero is the limit of the conditioned process
started from $r>0$ as $r\downarrow0$.

The modulus of variance-two three-dimensional Brownian motion
has generator $f''+2f'/r$. Adding vector drift $ae_1$ changes its
radial path law from entrance zero by the endpoint factor
\[
 e^{-a^2t/4}\frac{\sinh(ar/2)}{ar/2}.
\]
Indeed, the Brownian Girsanov density is
$e^{a B_t^{(1)}/2-a^2t/4}$; condition on the radial path and average
over its uniformly distributed endpoint direction. The resulting
radial generator is again $f''+a\coth(ar/2)f'$. This identifies
the process with the conditioned Brownian motion above \cite[Theorem~3]{RogersPitman1981}; see also
\cite{Alberts}. Set $a=a_\gamma$ and multiply by $\gamma$ to obtain
the stated representation. Independence follows from the wedge
construction.
\end{proof}

\subsection{Enhanced field data and embeddings}\label{subsec:enhanced-data}

Fix $s>0$ and $0<\sigma<1/2$ throughout. The field space is
$H^{-s}_{\mathrm{loc}}(\mathcal S)$. A distribution $\phi$ belongs
to this space if
\[
\xi\phi\in H^{-s}(\mathbb R^2)
\qquad\text{for every }\xi\in C_c^\infty(\mathcal S),
\]
where the localized distribution is extended by zero outside the
strip. Its topology is defined by
\[
\phi_n\to\phi
\quad\Longleftrightarrow\quad
\|\xi(\phi_n-\phi)\|_{H^{-s}(\mathbb R^2)}\to0
\quad\text{for every }\xi\in C_c^\infty(\mathcal S).
\]
The local space accommodates the linear growth of wedge fields
along the infinite strip.

For a wedge field $h_\gamma$, define its enhanced datum by
\[
Z=(\phi,\nu_0,\nu_\pi,\ell_0,\ell_\pi).
\]
The entries are as follows.
\begin{enumerate}
    \item
    $\phi=\gamma h_\gamma$ is the scaled field, regarded as an
    element of $H^{-s}_{\mathrm{loc}}(\mathcal S)$.

    \item
    $\nu_0,\nu_\pi$ are the GMC boundary length measures on
    $\mathbb R\times\{0\}$ and $\mathbb R\times\{\pi\}$,
    parametrized by their horizontal coordinate. We record them
    in $H^{-\sigma}_{\mathrm{loc}}(\mathbb R)$ with its strong
    local Sobolev topology. We will prove that the measures belong to
    this space almost surely for all sufficiently small $\kappa$.

    \item
    For $j\in\{0,\pi\}$,
    \[
    \ell_j(t)=\nu_j([t,\infty))
    \]
    is the cumulative boundary length measured from the
    $+\infty$ end. Wedge boundary lengths are finite towards
    $+\infty$ and infinite towards $-\infty$. Together with
    atomlessness and full support, this makes $\ell_j$ a
    continuous strictly decreasing bijection
    $\mathbb R\to(0,\infty)$. We record $\ell_j$ in
    $C_{\mathrm{loc}}(\mathbb R)$ with local uniform convergence.
\end{enumerate}
Thus
\[
Z\in\mathcal E:=
H^{-s}_{\mathrm{loc}}(\mathcal S)\times\bigl(H^{-\sigma}_{\mathrm{loc}}(\mathbb R)\bigr)^2
\times C_{\mathrm{loc}}(\mathbb R)^2,
\]
equipped with the product topology.

The measures are constructed by GMC at positive $\kappa$; they
are not pointwise exponentials of a distribution. At finite rate,
the limiting field has local $H^1$ regularity, and the enhanced
coordinates satisfy the classical identities
\[
\nu_j(dt)=e^{\phi(t,j)/2}\,dt,\qquad
\ell_j(t)=\int_t^\infty e^{\phi(u,j)/2}\,du,
\]
where boundary values are Sobolev traces.

We fix horizontal translations using one boundary length mark on
each surface:
\begin{equation}\label{eq:embeddingmarks}
\ell_{\mathrm{amb},0}(0)=1,\qquad
\ell_{1,0}(0)=1,\qquad
\ell_{2,\pi}(0)=1.
\end{equation}
These normalize the outer boundaries. The boundaries to be welded
are $\pi$ of component $1$ and $0$ of component $2$.
The upper and lower boundary measures of a given field are
generally different; normalizing one does not normalize the other.
We derive the field LDP first in the standard wedge embedding.
The boundary-distance embeddings in \eqref{eq:embeddingmarks} are obtained afterward
by contraction of the enhanced data, as proved in
Section~\ref{sec:enhanced}. The random translation is not assumed
negligible and does not amount to conditioning a boundary mass.

\subsection{Cutting maps and the joint data}\label{subsec:ABM}

Write $\mathcal K_\chi(\overline{\mathcal S}^{\,*})$ for the
space of nonempty compact subsets of the compactified strip,
equipped with the same Hausdorff metric induced by $d_\chi$.
The map $K\mapsto\overline\chi(K)$ identifies it isometrically
with $\mathcal K(\overline{\mathbb D})$. This larger space is
compact and allows subsequential limits that are not initially
known to be simple curves.

Write $\Phi$ for the scaled ambient field and $\phi_i$ for the
scaled component fields, with enhanced data $Z_{\mathrm{amb}}$
and $Z_i$. The curve $\eta$ separates $\mathcal S$ into $D_1$,
adjacent to $\theta=0$, and $D_2$, adjacent to $\theta=\pi$.
The cutting maps
\[
 F_i:\mathcal S\longrightarrow D_i\subset\mathcal S
\]
are conformal bijections preserving the two ends, with source and
target translations fixed by \eqref{eq:embeddingmarks}. At positive
$\kappa$, they obey the change-of-coordinate identity
\begin{equation}\label{eq:quantumcut}
 \phi_i=\Phi\circ F_i+(2+\kappa/2)\log|F_i'|.
\end{equation}
Their seam boundaries match according to length:
\begin{equation}\label{eq:weldexact}
 F_1\bigl(\ell_{1,\pi}^{-1}(l)+i\pi\bigr)
 =F_2\bigl(\ell_{2,0}^{-1}(l)\bigr),\qquad l>0.
\end{equation}
Boundary values are understood through the continuous extensions
associated with these simple-curve domains.

The geometric maps $G_i$ in \eqref{eq:geomark} depend only on the
curve, whereas the source normalization of $F_i$ uses the fields.
Outer boundary length correspondence gives the explicit relation
\begin{equation}\label{eq:derivedG}
 \begin{gathered}
 G_i(z)=F_i(z+s_i),\qquad
 s_i=\ell_{i,j_i}^{-1}\bigl(\ell_{\mathrm{amb},j_i}(0)\bigr),
 \qquad j_1=0,\quad j_2=\pi,\\
 s_1=0,\qquad G_1=F_1.
 \end{gathered}
\end{equation}
In general $s_2\ne0$, since the ambient upper length at zero need
not equal one. We therefore keep the distinction between $F_i$
and $G_i$, but record only $F_1,F_2$: the geometric maps are derived
from the enhanced data by \eqref{eq:derivedG}.

For maps
$H_n,H:\mathcal S\to\overline{\mathcal S}^{\,*}$, define locally
uniform convergence in $d_\chi$ by
\begin{equation}\label{eq:stripmaptopology}
 \sup_{z\in K}d_\chi(H_n(z),H(z))\longrightarrow0
 \quad\text{for every compact }K\subset\mathcal S.
\end{equation}
Equivalently, $\chi\circ H_n\circ\chi^{-1}$ converges locally
uniformly on $\mathbb D$ to $\chi\circ H\circ\chi^{-1}$, with
$\overline\chi$ understood at boundary values. Thus this is the
topology obtained by expressing the maps in disk coordinates. If the
limit is a holomorphic map into $\mathcal S$, it is also equivalent
to ordinary Euclidean local uniform convergence of the strip maps:
on each compact source set the limiting disk image stays inside
$\mathbb D$, where $\chi^{-1}$ is uniformly continuous locally.

Let $\mathcal F$ be the closure of the injective holomorphic maps
$\mathcal S\to\mathcal S$ in
$C_{\mathrm{loc}}(\mathcal S;\overline{\mathcal S}^{\,*})$ with
this topology. It is compact by the normal-family theorem applied
to their disk representatives. Its elements include constant maps
with values on the compactified target boundary, so degeneration
is allowed before finite-rate stability is applied.

The joint data and state space are
\[
 A=Z_{\mathrm{amb}}\in\mathcal E,\qquad
 B=(Z_1,Z_2)\in\mathcal E^2,\qquad
 M=(F_1,F_2,\eta)\in\mathcal M,
\]
where
\[
 \mathcal M=\mathcal F^2
       \times\mathcal K_\chi(\overline{\mathcal S}^{\,*}).
\]
The product topology combines \eqref{eq:stripmaptopology} for
the two maps with $d_{\mathrm H}^{\chi}$ for the curve. All
coordinates are therefore expressed on the strip; $\chi$ is
used to define their topology and prove compactness. Bounded disk
representatives will be introduced locally inside proofs when
area or boundary trace estimates are needed. They are not
additional coordinates of the datum.
For later use, denote the enhanced wedge rate by
$\mathbf I_w^{(j)}(Z)$. It equals $I_w^{(j)}(\phi)$ when the
field has finite energy and the remaining coordinates are its
classical measures and cumulative lengths; it is infinite
otherwise. The enhanced marginal LDPs will have rates
\[
 I(A)=\mathbf I_W^{(0)}(Z_{\mathrm{amb}}),\qquad
 K(B)=\mathbf I_{W_1}^{(0)}(Z_1)
           +\mathbf I_{W_2}^{(\pi)}(Z_2).
\]
The target joint rate is
\begin{equation}\label{eq:ABMtarget}
 L(A,B,M)=
 \begin{cases}
 K(B),&\begin{gathered}
 \text{at the unique normalized welding configuration}\\
 \text{above a }B\text{ of finite rate},
 \end{gathered}\\[4pt]
 +\infty,&\text{otherwise}.
 \end{cases}
\end{equation}
Section~\ref{sec:curve} establishes this joint LDP and then projects
to the interface. The compact enlargement of the geometric space
is used for subsequential LDPs, not as a claim that welding is
continuous on all of $\mathcal E^2$.
\subsection{The field large deviation principle}

For a finite-energy profile $\phi$, let $p$ be its angular mean
and write $v=\phi-p$. The energy in \eqref{eq:energy} is
\[
 \frac14\int_\R|p'+w|^2\dd t
 +\frac1{4\pi}\int_{\Sstrip}|\nabla v|^2\dd^2z.
\]
Here $v\in H^1(\Sstrip)$ by the angular Poincare inequality. 
\begin{proposition}[Standard-embedding field LDP]\label{prop:raw}
The standard-embedding field $\phi_\kappa$ satisfies a good LDP in
$H^{-s}_{\mathrm{loc}}(\mathcal S)$ with rate
\begin{equation}\label{eq:rawrate}
 R_w(\phi)=\frac14\int_\R|r'+w|^2\dd t
 +\frac1{4\pi}\int_{\Sstrip}|\nabla v|^2\dd^2z,
 \qquad \phi=r+v.
\end{equation}
Its domain has $r\in AC_{\rm loc}$, $r(0)=0$,
$r(t)\ge0$ for $t\le0$, $r'+w\in L^2(\R)$, and
$v\in H^1(\Sstrip)$ with zero angular average. It is infinite otherwise.
\end{proposition}
\begin{proof}
For the Gaussian construction we may first use the stronger auxiliary
space $C_{\mathrm{loc}}(\mathbb R;H_N^{-s}(0,\pi))$.
Here $H_N^{-s}(0,\pi)$ is the spectral Sobolev space defined by
the Neumann Laplacian on $(0,\pi)$, using the constant mode and
the cosine modes.
Its elements are continuous distribution-valued paths in the
horizontal coordinate. There is a continuous map
\[
 C_{\mathrm{loc}}(\mathbb R;H_N^{-s}(0,\pi))\longrightarrow H^{-s}_{\mathrm{loc}}(\mathcal S)
\]
obtained by viewing such a path as a distribution on the strip.
Continuity follows by localization and the continuous embedding of
$L^2$ in the horizontal variable with values in angular $H_N^{-s}$
into two-dimensional local $H^{-s}$. We prove the LDP in this
auxiliary space and then contract to $H^{-s}_{\mathrm{loc}}(\mathcal S)$.
Apply the Gaussian Schilder theorem on the local path and field spaces
to the independent Brownian motions and $X^\perp$. The deterministic
drift correction is at most $\kappa T/2$ on each time interval of
length $T$, so it is exponentially equivalent in the local uniform
topology. Lemma~\ref{lem:radial} and contraction give the negative half.
For completeness, the contraction cost of a nonnegative path $y$, $y(0)=0$,
on $[0,T]$ is
\[
 \inf_{|x|=y}\frac14\int_0^T|x'-we_1|^2
 =\frac14\int_0^T(y'-w)^2.
\]
Indeed $|x'|\ge|y'|$ a.e. and $x_1(T)\le y(T)$ give the lower bound
after expansion; $x=ye_1$ attains it. Independent Gaussian energies
add, and the increasing finite-interval energies give
\eqref{eq:rawrate}. Goodness follows from the Gaussian contraction
theorem in the auxiliary space. The lateral GFF is a Radon Gaussian
on $C_{\mathrm{loc}}(\mathbb R;H_N^{-s}(0,\pi))$; its nonconstant
angular modes are stationary Ornstein--Uhlenbeck processes with
variances proportional to $1/n$. Their weighted sums give continuous
$H_N^{-s}$ paths for every $s>0$.
Finally contract to $H^{-s}_{\mathrm{loc}}(\mathcal S)$. On the finite-rate domain, the resulting
distribution determines $r$ as its angular mean and $v=\phi-r$
uniquely. The contraction therefore leaves both the displayed energy
and its admissibility conditions unchanged. Goodness is preserved
by this continuous contraction.
\end{proof}

\section{Boundary measures and LDP for $B$}\label{sec:enhanced}

We prove the enhanced LDP in the space $\mathcal E$ of
Subsection~\ref{subsec:enhanced-data}, and hence the marginal LDPs
for $A$ and $B$. We first work in the standard wedge embedding;
the final step applies the length normalization
\eqref{eq:embeddingmarks} by contraction.

\begin{theorem}[Enhanced wedge LDP]\label{thm:wedge}
The standard-embedding enhanced datum $Z_\kappa$ has
a good LDP in $\mathcal E$ with rate \eqref{eq:rawrate} on the classical boundary measure
and cumulative-length graph, and infinite rate elsewhere.

For $j=0$ or $\pi$, let $t_\kappa$ solve
$\ell_{\kappa,j}(t_\kappa)=1$, and translate every coordinate by
$t_\kappa$. This unit boundary-distance embedding has a good enhanced
LDP with rate $I_w^{(j)}$ in \eqref{eq:energy} on the same classical
graph, and infinity elsewhere. Its unique zero is
\eqref{eq:zero} with its classical measures and cumulative lengths.
\end{theorem}

The proof combines local chaos approximation with a tail estimate at
the finite-length end. The local estimates are proved first; the
normalization and the cumulative coordinates are assembled at the
end of the section.

\subsection{Local chaos approximation}

\begin{lemma}[Hilbert-valued boundary GMC estimate]\label{lem:HC}
For a smooth compactly supported boundary cutoff $\zeta\ge0$ and
$0<\beta^2<4\sigma$, uniformly for $0<\kappa\le\beta^2$,
\[
 \E\|\zeta\nu_{\beta^2,j}^\perp\|_{H^{-\sigma}}^2<\infty,
 \qquad
 \left\|\|\zeta\nu_{\kappa,j}^\perp\|_{H^{-\sigma}}\right\|_
       {L^{1+\beta^2/\kappa}}\le C_{\beta,\zeta}.
\]
If the lateral chaos is replaced by its conditional expectation on
the first $m$ smooth Gaussian coordinates, the corresponding difference
has the same high-moment bound with a constant tending to zero as
$m\to\infty$.
\end{lemma}
\begin{proof}
We first use full Wick subtraction. Let $M_a$ denote the cutoff
lateral boundary chaos, valued in $H^{-\sigma}(\mathbb R)$, with parameter
$a$, whose exponential coefficient is $a/2$. Thus
$M_{\sqrt\kappa}=\zeta\nu_{\kappa,j}^\perp$ in this convention.

\emph{The fixed-parameter second moment.}
Let $K_\sigma$ be the kernel of $(1-\Delta)^{-\sigma}$ on
$\R$. For a positive cutoff measure $M$,
\[
 \|M\|_{H^{-\sigma}(\mathbb R)}^2
 =\iint K_\sigma(t-u)\,M(\dd t)M(\dd u).
\]
Near zero the Bessel kernel $K_\sigma$ is comparable to
$|t-u|^{2\sigma-1}$. The boundary covariance
$-2\log|t-u|+O(1)$ and the Gaussian exponential identity give a
two-point GMC factor bounded by $C|t-u|^{-\beta^2/2}$ on the
cutoff support. Consequently
\[
 \E\|M_\beta\|_{H^{-\sigma}(\mathbb R)}^2
 \le C\iint\zeta(t)\zeta(u)
       |t-u|^{2\sigma-1-\beta^2/2}\dd t\dd u+C<\infty,
\]
because $2\sigma-\beta^2/2>0$. These calculations can first be made
for smooth regularizations of the field. The same integrable bounds
on their mixed second moments give convergence in $L^2(H^{-\sigma}(\mathbb R))$
and justify the identities for the limiting chaos.

\emph{Mehler's formula and the high-moment estimate.}
Represent the lateral field by independent standard Gaussian
coordinates $\xi=(\xi_n)_{n\ge1}$. For $0\le r\le1$ define
the Gaussian noise operator $T_r$ by
\begin{equation}\label{eq:Mehler}
 (T_rF)(\xi)=\E_{\xi'}\left[
 F\bigl(r\xi+\sqrt{1-r^2}\,\xi'\bigr)\right],
\end{equation}
where $\xi'$ is an independent copy and only $\xi'$ is averaged.
This is Mehler's formula for the Ornstein--Uhlenbeck semigroup,
with $r=e^{-t}$ at semigroup time $t$. It replaces part of the
original Gaussian field by independent noise while preserving its
law. For $H^{-\sigma}(\mathbb R)$-valued $F$ the expectation is a Bochner
expectation, so the same definition applies to our random measures.

The key consequence is
\begin{equation}\label{eq:Mehlerchaos}
 T_rM_\beta=M_{r\beta}.
\end{equation}
Indeed, for a centered Gaussian regularization $Y$ of the field,
an independent copy $Y'$, and any real $a$, Gaussian integration gives
\[
 \E_{Y'}\exp\left\{
 a(rY+\sqrt{1-r^2}Y')-\frac{a^2}{2}\E Y^2\right\}
 =\exp\left\{arY-\frac{a^2r^2}{2}\E Y^2\right\}.
\]
Use $a=\beta/2$,
multiply by the cutoff, and pass to the $L^2(H^{-\sigma}(\mathbb R))$ limit.
This proves \eqref{eq:Mehlerchaos}; no exponential of an
unregularized distribution is used in this calculation.

Scalar Gaussian hypercontractivity \cite[Chapter 5]{Janson} states that
\[
 \|T_rf\|_{L^q}\le\|f\|_{L^p},\qquad
 1<p\le q<\infty,\qquad r^2\le\frac{p-1}{q-1}.
\]
Jensen's inequality gives
$\|T_rF\|_{H^{-\sigma}(\mathbb R)}\le T_r(\|F\|_{H^{-\sigma}(\mathbb R)})$ pointwise.
Applying the scalar estimate to $\|F\|_{H^{-\sigma}(\mathbb R)}$ therefore yields
the same bound for the $L^q$ norm of the Hilbert norm. Set
\[
 p=2,\qquad r=\sqrt\kappa/\beta,\qquad q=1+\beta^2/\kappa.
\]
Then $r\le1$, $q\ge2$, and $r^2=(p-1)/(q-1)$. Hence
\begin{equation}\label{eq:HChypercontractive}
 \left\|\|M_{\sqrt\kappa}\|_{H^{-\sigma}(\mathbb R)}\right\|_{L^q}
 \le\left\|\|M_\beta\|_{H^{-\sigma}(\mathbb R)}\right\|_{L^2}<\infty.
\end{equation}
The right side is independent of $\kappa$, proving the claimed
high-moment bound in the Wick convention.

\emph{Finite-coordinate approximation.}
Choose the Gaussian coordinates with smooth Cameron--Martin modes,
put $\mathcal F_m=\sigma(\xi_1,\ldots,\xi_m)$, and write
$P_mF=\E[F\mid\mathcal F_m]$. Independence of the coordinates
in \eqref{eq:Mehler} shows that $P_mT_r=T_rP_m$: both compositions
average over the unused coordinates and the independent added noise.
Thus
\[
 M_{\sqrt\kappa}-P_mM_{\sqrt\kappa}
 =T_r(M_\beta-P_mM_\beta).
\]
Applying the Hilbert-valued estimate above to this signed difference
gives the explicit uniform error bound
\begin{equation}\label{eq:HCapprox}
 \left\|\|M_{\sqrt\kappa}-P_mM_{\sqrt\kappa}\|_{H^{-\sigma}(\mathbb R)}
       \right\|_{L^{1+\beta^2/\kappa}}
 \le\left\|\|M_\beta-P_mM_\beta\|_{H^{-\sigma}(\mathbb R)}\right\|_{L^2}
 \longrightarrow0.
\end{equation}
The convergence is Hilbert-valued martingale convergence:
$M_\beta\in L^2(H^{-\sigma}(\mathbb R))$ and the coordinate sigma-fields
increase to the sigma-field of the lateral field. In particular,
the error tends to zero uniformly for $0<\kappa\le\beta^2$.

Finally return to geometric normalization. On the cutoff support
it multiplies the Wick chaos by a deterministic smooth function
$d_a$, with $d_a\to1$ in smooth norm as $a\to0$.
After localization by a fixed larger cutoff, multiplication by $d_a$
is bounded on $H^{-\sigma}(\mathbb R)$, uniformly for $0\le a\le\beta$.
It commutes with $P_m$ because it is deterministic. Applying this
multiplier to the chaos and to its approximation preserves both
\eqref{eq:HChypercontractive} and \eqref{eq:HCapprox}, up to a
fixed constant. This proves the lemma in the normalization used here.
\end{proof}

\begin{proposition}[Local exponential approximation]\label{prop:local}
As $\kappa\downarrow0$, the standard-embedding random tuple
\[
 (\phi_\kappa,\nu_{\kappa,0},\nu_{\kappa,\pi})
 \in H^{-s}_{\mathrm{loc}}(\mathcal S)\times\bigl(H^{-\sigma}_{\mathrm{loc}}(\mathbb R)\bigr)^2
\]
satisfies a good LDP at speed $\kappa^{-1}$ with rate
\begin{equation}\label{eq:localjointrate}
 \mathcal I_w(\phi,m_0,m_\pi)
 =\begin{cases}
 R_w(\phi),
       &\begin{gathered}
          \phi=r+v\text{ is admissible},\\
          m_j=e^{\phi(t,j)/2}\dd t\quad(j=0,\pi),
        \end{gathered}\\[6pt]
 +\infty,&\text{otherwise}.
 \end{cases}
\end{equation}
Admissibility is as in Proposition~\ref{prop:raw}. The boundary
values of an admissible $\phi$ are its local Sobolev traces.
The classical measures in \eqref{eq:localjointrate} belong to
the stated spaces. The random measures belong to these spaces
almost surely for all sufficiently small $\kappa$.
\end{proposition}
\begin{proof}
\emph{Local spaces and almost sure membership.}
Choose nonnegative cutoffs $\zeta_n$ on $\R$,
equal to one on successively larger compact intervals and
with each support contained in the next region where the cutoff
is one. These intervals exhaust $\R$. A metric on $H^{-\sigma}_{\mathrm{loc}}(\mathbb R)$ is
\[
 d_{-\sigma,\mathrm{loc}}(m,\widetilde m)
 =\sum_{n\ge1}2^{-n}
   \bigl(1\wedge\|\zeta_n(m-\widetilde m)\|_{H^{-\sigma}}\bigr).
\]
Multiplication by smooth cutoffs is continuous on Sobolev spaces,
so this metric gives the topology fixed in
Subsection~\ref{subsec:enhanced-data}.
The compatible cutoff distributions form a closed subspace of a
countable product of separable Hilbert spaces, proving the Polish
assertion. Positivity is also closed under this convergence; one
could equivalently use the closed subspaces of positive Radon
measures that belong to these local Sobolev spaces.

Fix $\sigma<\alpha<1/2$ and write $r_\kappa=\gamma R_\gamma$.
On compact intervals $r_\kappa\in C^\alpha$ almost surely.
Lemma~\ref{lem:HC}, with any fixed $0<\beta^2<4\sigma$ and
$\kappa\le\beta^2$, puts each cutoff lateral measure in
$H^{-\sigma}$ almost surely. Multiplication by a $C^\alpha$
function is bounded on $H^{-\sigma}$, by duality from the
corresponding $H^\sigma$ multiplier estimate. Thus multiplication
by $e^{r_\kappa/2}$ gives the asserted local membership of the full
boundary measures, simultaneously over the countable exhaustion.

\emph{Exponentially good finite-coordinate approximations.}
Use the coordinates of Lemma~\ref{lem:HC} and put
$X_m=\E[X^\perp\mid\mathcal F_m]$. Write
$K_m(t,j)=\E[X_m(t,j)^2]$. With Wick normalization the lateral approximations have
densities
\[
 \exp\{\sqrt\kappa X_m(t,j)/2-\kappa K_m(t,j)/8\}.
\]
They are precisely the conditional expectations of the lateral
boundary chaoses: the unobserved Gaussian modes average to one after
Wick subtraction. Multiply these densities by $e^{r_\kappa/2}$
and by the deterministic geometric normalization factors. Call the
resulting full measures $\nu^{(m)}_{\kappa,j}$; retain the original
field $\phi_\kappa$ in the approximating tuple.

For every compact interval $K$ and every $A<\infty$ there is
$L=L(A,K)$ such that
\begin{equation}\label{eq:radialHoldertight}
 \limsup_{\kappa\downarrow0}\kappa\log
 \Pp\{\|r_\kappa\|_{C^\alpha(K)}>L\}\le-A.
\end{equation}
Indeed, Lemma~\ref{lem:radial} represents the negative half as
the norm of a Gaussian vector path with deterministic drift; the
norm map bounds its H\"older seminorm by that of the vector path.
Gaussian tails for Brownian $C^\alpha$ norms give the estimate on
both halves. On the event in which this radial norm is at most $L$,
the multiplier estimate bounds a cutoff full-measure error by
$C_L$ times the cutoff lateral error. If
$q_\kappa=1+\beta^2/\kappa$ and $\varepsilon_m\to0$ is the
$L^{q_\kappa}(H^{-\sigma})$ error from Lemma~\ref{lem:HC},
Markov's inequality gives, for a boundary cutoff $\zeta$,
\[
 \begin{split}
 &\Pp\{\|\zeta(\nu_{\kappa,j}
                  -\nu^{(m)}_{\kappa,j})\|_{H^{-\sigma}}>\delta\}\\
 &\quad\le\Pp\{\|r_\kappa\|_{C^\alpha(K)}>L\}
              +(C_L\varepsilon_m/\delta)^{q_\kappa},
 \end{split}
\]
where $K$ contains the cutoff support. Therefore the logarithmic
upper bound is at most
$\max\{-A,\beta^2\log(C_L\varepsilon_m/\delta)\}$.
First send $m\to\infty$ with $L$ fixed, then send $A\to\infty$.
A finite union bound over the two boundaries and the summable tails
in $d_{-\sigma,\mathrm{loc}}$
prove, for the product metric $d$ on $H^{-s}_{\mathrm{loc}}(\mathcal S)\times\bigl(H^{-\sigma}_{\mathrm{loc}}(\mathbb R)\bigr)^2$,
\begin{equation}\label{eq:localexpapprox}
 \lim_{m\to\infty}\limsup_{\kappa\downarrow0}\kappa\log
 \Pp\{d((\phi_\kappa,\nu_{\kappa,0},\nu_{\kappa,\pi}),
 (\phi_\kappa,\nu^{(m)}_{\kappa,0},\nu^{(m)}_{\kappa,\pi}))
                         >\delta\}=-\infty.
\end{equation}

\emph{Finite-coordinate contraction and the limiting graph.}
The Gaussian Schilder LDP used in Proposition~\ref{prop:raw}
can be applied jointly to the underlying Brownian paths, the
lateral field, and any finite collection of its Gaussian coordinates.
One may choose these coordinates from a generating family of smooth
test-function pairings by Gram--Schmidt, so that their scaled values
are continuous linear observations on the local field space.
For fixed $m$ the displayed densities are continuous functions of
these finitely many scaled coordinates and of the radial path.
The variance corrections tend to zero smoothly on each cutoff;
the same is true of the deterministic normalization corrections.
Contraction, or equivalently its version for such converging maps,
therefore gives an LDP for each approximating tuple.

On a finite-rate profile $\phi=r+v$, let $v_m$ be the
Cameron--Martin orthogonal projection of $v$ onto the first $m$
lateral modes. The deterministic approximating measures are
\[
 e^{(r+v_m)(t,j)/2}\dd t,\qquad j\in\{0,\pi\}.
\]
The underlying Schilder cost contracts to \eqref{eq:rawrate}:
the radial contraction is exactly that computed in
Proposition~\ref{prop:raw}, and the added observations have no
independent cost. We now verify convergence of these deterministic
measures uniformly on bounded energy sets, as required to identify
the limiting rate.

If the energy is at most $A$, then $r(0)=0$ and the bound on
$r'+w$ give uniform local $H^1$ bounds for $r$. The zero angular
average and Poincare's inequality give a global $H^1$ bound for
$v$ and for all $v_m$. The projections converge strongly in the
Cameron--Martin space for each $v$. Uniformly over its bounded
balls, $v_m-v$ tends to zero in boundary $L^2$ trace norm:
the localized trace map from $H^1$ to boundary $L^2$ is compact,
and orthogonal projections
converging strongly to the identity converge uniformly after any
such compact map. The exponential trace inequality gives,
on each boundary cutoff, uniform $L^p$ bounds for the exponential
densities for every finite $p$ \cite{Ozawa1995}. Hence convergence in measure and
uniform integrability give convergence of the densities in local
$L^2$, uniformly over the energy sublevel. Since
$\|f\|_{H^{-\sigma}}\le\|f\|_{L^2}$, this proves their convergence
in the stated strong local Sobolev topologies. It also proves that
the limiting classical measures belong to those spaces.

The same compactness argument applies to sequences of admissible
profiles of bounded energy converging in $H^{-s}_{\mathrm{loc}}(\mathcal S)$: their local weak
$H^1$ limits and strong $L^2$ limits, including boundary traces,
are determined by the limiting field. Thus the classical-measure
graph map is continuous on every rate sublevel. The exponentially
good approximation theorem applied to \eqref{eq:localexpapprox}
now gives exactly \eqref{eq:localjointrate}. Its sublevels are
the continuous graph images of the compact sublevels of
Proposition~\ref{prop:raw}, so the rate is good.

For completeness, exponential tightness of the local measure
coordinates can also be checked directly. Fix $0<\sigma_0<\sigma$.
Positivity of the Bessel kernel gives
\[
 \|\zeta\nu_{\kappa,j}\|_{H^{-\sigma_0}}
 \le e^{\frac12\sup_{\supp\zeta}\gamma R_\gamma}
       \|\zeta\nu_{\kappa,j}^\perp\|_{H^{-\sigma_0}}.
\]
Independence, radial Gaussian supremum bounds, and
Lemma~\ref{lem:HC} imply, for some $c>0$,
\begin{equation}\label{eq:normmoment}
 \limsup_{\kappa\to0}\kappa\log
 \E\|\zeta\nu_{\kappa,j}\|_{H^{-\sigma_0}}^{c/\kappa}<\infty.
\end{equation}
Norm balls in $H^{-\sigma_0}$ with fixed compact support have compact
closure in $H^{-\sigma}$. Equation~\eqref{eq:normmoment} and
Markov's inequality therefore give exponential tightness on each
cutoff. The moment estimates used here in fact give
$\E\|\zeta_n\nu_{\kappa,j}\|_{H^{-\sigma_0}}^{c/\kappa}
\le\exp(C_n/\kappa)$, with a common $c>0$ and common
small-$\kappa$ range; only $C_n$ depends on the cutoff.
For any target rate $A$, choose bounds $R_n$ with
$c\log R_n\ge C_n+A+n$. The exceptional probability on the
$n$th cutoff is then at most $\exp[-(A+n)/\kappa]$.
A union bound and diagonal compactness give exponential tightness
in the countable local topology.
\end{proof}

\begin{remark}
The error of two chaos approximations is signed. Positivity of the
Bessel kernel justifies the norm bound for a positive measure in
\eqref{eq:normmoment}; it does not justify multiplying a signed error
by a merely bounded radial function. The $C^\alpha$ multiplier bound
in the proof is needed for that error.
\end{remark}

\subsection{Finite-end tails and length coordinates}

\begin{lemma}[Exponential boundary-tail bound]\label{lem:tail}
For every $w>0$ there is $c>0$ such that, for all sufficiently small
$\kappa$ and $T\ge1$,
\[
 \|\nu_{\kappa,j}([T,\infty))\|_{L^{c/\kappa}}
 \le C e^{-wT/4}.
\]
In particular these tails are exponentially negligible at speed
$\kappa^{-1}$ as $T\to\infty$.
\end{lemma}
\begin{proof}
Choose a fixed $\beta$ with $\beta^2<2$ and then
$0<c<\min(\beta^2,w)/2$. Put $q=c/\kappa\ge1$.
Stationarity and Lemma~\ref{lem:HC}, or its scalar total-mass version,
bound the lateral unit-cell $L^q$ norm uniformly. The radial value at
the integer $n$ contributes
\[
 \|e^{\gamma R_\gamma(n)/2}\|_{L^q}
 =\exp\left[-\left(\frac w2-\frac{\kappa}{4}
                                 -\frac{q\kappa}{4}\right)n\right].
\]
The independent Brownian increments on a unit cell have a uniformly
bounded exponential $L^q$ norm since $q\kappa=c$. Thus
\[
 \|\nu_{\kappa,j}([n,n+1])\|_{L^q}
 \le C\exp[-(w/2-(c+\kappa)/4)n].
\]
For small $\kappa$ the exponent is at least $w/4$. Sum the bounds
using Minkowski. Markov gives
\[
 \limsup_{\kappa\to0}\kappa\log
 \Pp\{\nu_{\kappa,j}([T,\infty))>\varepsilon\}
 \le c\log C-cwT/4-c\log\varepsilon,
\]
which tends to $-\infty$.
\end{proof}

\begin{lemma}[Classical cumulative lengths]\label{lem:classcdf}
If $\phi=-wt+u$ and $\int|\nabla u|^2<\infty$, then both classical
boundary measures have finite mass on $[T,\infty)$ and infinite mass
on $(-\infty,T]$. Their cumulative functions are continuous strictly
decreasing bijections $\R\to(0,\infty)$.
\end{lemma}
\begin{proof}
Write $u=p+v$ with zero angular average for $v$. Then $p'\in L^2$
and $v\in H^1$ by the angular Poincare inequality. Hence
$p(t)=o(|t|)$ at both ends. The trace exponential inequality on unit
rectangles and Jensen's inequality give upper and lower bounds for
the unit-cell exponential trace integrals whose logarithms are
$o(|n|)$; the global $H^1$ bound of $v$ and the $L^2$ bound on $p'$
make these bounds uniform over cells, apart from $p(n)$.
Thus the logarithm of the boundary mass of $[n,n+1]$ is
$-wn/2+o(|n|)$. Summation proves both assertions. Densities are positive
a.e. and locally integrable, proving continuity and strictness.
\end{proof}

\subsection{Assembly and choice of embedding}

\begin{proof}[Proof of Theorem~\ref{thm:wedge}]
Put
\[
 L_j^\phi(t)=\int_t^\infty e^{\phi(u,j)/2}\dd u.
\]
The two cumulative coordinates carry local uniform convergence.
We first prove the LDP in $\mathcal E$, then apply the length-one
translation.

\emph{1. Existence and truncation of cumulative lengths.}
For sufficiently small $\kappa$, Lemma~\ref{lem:tail} and the
atomless full-support property of subcritical GMC give
\[
 0<\nu_{\kappa,j}([T,\infty))<\infty,\qquad
 \nu_{\kappa,j}(\{t\})=0,\qquad
 \nu_{\kappa,j}((a,b))>0\quad(a<b).
\]
At the other end, for fixed $\kappa$,
\[
 \frac{r_\kappa(-t)}t\longrightarrow w-\kappa/2>0
 \quad\text{a.s.},\qquad r_\kappa=\gamma R_\gamma.
\]
Stationarity and ergodicity of the lateral chaos imply, for some
$\delta>0$, infinitely many integers $n\to\infty$ with
$\nu^\perp_{\kappa,j}([-n-1,-n])\ge\delta$. On these cells the
radial multiplier grows exponentially, so
\[
 \nu_{\kappa,j}(( -\infty,T])=\infty,\qquad
 \ell_{\kappa,j}(-\infty)=\infty,\quad
 \ell_{\kappa,j}(+\infty)=0.
\]
Thus $\ell_{\kappa,j}$ is a continuous strictly decreasing
bijection $\R\to(0,\infty)$. We use here the usual GMC
nondegeneracy and support properties summarized in \cite{RV}.

For $T>R$ and $|t|\le R$, define
\[
 \ell^{T}_{\kappa,j}(t)=\nu_{\kappa,j}([t,T]).
\]
Choose $\zeta\in C_c^\infty(\R)$ equal to one on $[-R,T]$.
For $0<\sigma<1/2$, the indicator functions satisfy
\[
 \sup_{|t|\le R}\|\mathbf1_{[t,T]}\|_{H^\sigma}\le C_{R,T},
 \qquad
 \|\mathbf1_{[t,T]}-\mathbf1_{[u,T]}\|_{H^\sigma}
 \le C_{R,T}|t-u|^{1/2-\sigma}.
\]
Consequently the pairing with a local $H^{-\sigma}$ distribution
is a continuous function of $t$, and
\begin{equation}\label{eq:cdfSobolevcontinuity}
 \sup_{|t|\le R}|\ell^T_{n,j}(t)-\ell^T_j(t)|
 \le C_{R,T}\|\zeta(\nu_{n,j}-\nu_j)\|_{H^{-\sigma}}
 \longrightarrow0.
\end{equation}
For our positive atomless measures this pairing is the usual interval
mass. Hence Proposition~\ref{prop:local} and contraction give the
joint LDP with the truncated cumulative coordinates.

The truncation error is particularly simple:
\[
 \sup_{|t|\le R}|\ell_{\kappa,j}(t)
                    -\ell^T_{\kappa,j}(t)|
 =\nu_{\kappa,j}([T,\infty)).
\]
Lemma~\ref{lem:tail} therefore gives, for $\varepsilon>0$,
\begin{equation}\label{eq:cdfexptail}
 \limsup_{\kappa\downarrow0}\kappa\log
 \Pp\!\left\{\sup_{|t|\le R}|\ell_{\kappa,j}(t)
                 -\ell^T_{\kappa,j}(t)|>\varepsilon\right\}
 \le c\log(C/\varepsilon)-cwT/4
 \xrightarrow[T\to\infty]{}-\infty.
\end{equation}
To regard the approximations as elements of $C_{\rm loc}(\R)$,
set $\ell^T(t)=0$ for $t\ge T$. Continuity of this extension
also follows from the indicator estimate above.

\emph{2. The deterministic tail bound and the enhanced rate.}
The same truncation must converge uniformly on rate sublevels.
If $R_w(\phi)\le A$, write $\phi=r+v$. Then
\[
 r(0)=0,\qquad \|r'+w\|_{L^2}\le2\sqrt A,
 \qquad \|v\|_{H^1(\Sstrip)}\le C\sqrt A,
\]
and Cauchy--Schwarz gives
\[
 r(t)=-wt+\int_0^t(r'(u)+w)\dd u
       \le-wt+2\sqrt{At}\qquad(t\ge0).
\]
The exponential trace inequality on unit rectangles gives
$\int_n^{n+1}e^{v(t,j)/2}\dd t\le C_A$, uniformly in $n$ and
in the rate sublevel. Thus, for integers $n\ge0$,
\[
 \int_n^{n+1}e^{\phi(t,j)/2}\dd t
 \le C_A\exp\{-wn/2+\sqrt{A(n+1)}\}
 \le C_{A,w}e^{-wn/4}.
\]
Summing yields the explicit uniform estimate
\begin{equation}\label{eq:classicaluniformtail}
 \sup_{R_w(\phi)\le A}L_j^\phi(T)
 \le C_{A,w}e^{-wT/4}\longrightarrow0.
\end{equation}
Equations~\eqref{eq:cdfSobolevcontinuity} and
\eqref{eq:classicaluniformtail} show that adjoining $L_j^\phi$
is continuous on each compact rate sublevel of
Proposition~\ref{prop:local}. Together with
\eqref{eq:cdfexptail}, the exponentially good approximation theorem
therefore gives the enhanced rate
\[
 \mathcal R_w(\phi,m_0,m_\pi,\ell_0,\ell_\pi)
 =\begin{cases}
 R_w(\phi),&
 \begin{gathered}
 R_w(\phi)<\infty,\quad m_j=e^{\phi(t,j)/2}\dd t,\\
 \ell_j=L_j^\phi\quad(j=0,\pi),
 \end{gathered}\\[4pt]
 +\infty,&\text{otherwise}.
 \end{cases}
\]
Its sublevels are continuous graph images of compact sublevels,
so it is good. The local uniform metrics use summable weights
$2^{-n}$ on compact windows; for any fixed metric error, only
finitely many windows matter. Thus the finite-window estimates
above also give the asserted countable local enhancement.

\emph{3. Continuity of the length-one translation.}
For finite-rate data, Lemma~\ref{lem:classcdf} gives a unique
\[
 t_j=\ell_j^{-1}(1).
\]
Suppose $\ell_{n,j}\to\ell_j$ locally uniformly. For every
$\varepsilon>0$,
\[
 \ell_j(t_j-\varepsilon)>1>\ell_j(t_j+\varepsilon),
\]
so, for all large $n$, the same strict inequalities hold for
$\ell_{n,j}$. If these functions are decreasing, their roots obey
\begin{equation}\label{eq:lengthrootbracket}
 t_j-\varepsilon<t_{n,j}<t_j+\varepsilon,
 \qquad t_{n,j}\longrightarrow t_j.
\end{equation}
For contraction on the entire coordinate space, extend the root
map to arbitrary continuous functions by choosing the root of
smallest absolute value, breaking ties by choosing the smaller
root, and choosing zero if there is no root. Bracketing guarantees
a root near $t_j$. The chosen roots are therefore bounded; local
uniform convergence excludes any limit other than $t_j$.
This measurable extension is consequently continuous at every
finite-rate datum as well.

Define simultaneous translation of the coordinates by
\[
 (\tau_t\phi)(x,\theta)=\phi(x+t,\theta),\qquad
 (\tau_tm)(B)=m(B+t),\qquad
 (\tau_t\ell)(x)=\ell(x+t),
\]
where $B+t=\{x+t:x\in B\}$ for boundary sets $B\subset\mathbb R$.
Translation is jointly continuous in the local field, Sobolev,
and uniform topologies. The normalization map is therefore
\[
 \mathcal N_j(z)=\tau_{t_j(z)}z,\qquad
 (\tau_{t_j}\ell_j)(0)=1,
\]
continuous at every finite-rate point. The contraction principle
in this form gives a good LDP with rate
\begin{equation}\label{eq:normalizedcontraction}
 \widehat{\mathcal R}_w^{(j)}(\widehat z)
 =\inf\{\mathcal R_w(z):\mathcal N_j(z)=\widehat z\}.
\end{equation}

\emph{4. Evaluation of the contraction infimum.}
For any finite-energy profile define the translation-invariant energy
\[
 \mathcal E_w(\phi)=\frac1{4\pi}
       \int_{\Sstrip}|\nabla(\phi+wt)|^2\dd^2z,
 \qquad \mathcal E_w(\tau_a\phi)=\mathcal E_w(\phi).
\]
On the raw admissible domain, $\mathcal E_w=R_w$. Every preimage
in \eqref{eq:normalizedcontraction} thus has the same energy as
its normalized field, and its normalized cumulative coordinate
satisfies $L_j^{\widehat\phi}(0)=1$. This proves the lower bound
by $I_w^{(j)}(\widehat\phi)$.

Conversely, suppose $\mathcal E_w(\widehat\phi)<\infty$ and
$L_j^{\widehat\phi}(0)=1$. Its angular mean $\widehat r$ satisfies
\[
 \widehat r(t)=-wt+o(|t|)\quad(t\to\pm\infty).
\]
Let $a$ be its first zero. Continuity and the limits at both ends give
\[
 \widehat r(a)=0,\qquad\widehat r(t)\ge0\quad(t\le a).
\]
The translated profile $\phi=\tau_a\widehat\phi$ is raw admissible,
with $R_w(\phi)=\mathcal E_w(\widehat\phi)$. Moreover
\[
 L_j^\phi(x)=L_j^{\widehat\phi}(x+a),\qquad
 (L_j^\phi)^{-1}(1)=-a,\qquad
 \tau_{-a}\phi=\widehat\phi.
\]
It is therefore a preimage attaining the infimum. Hence the normalized
rate is $I_w^{(j)}$ on its classical enhanced graph and infinity
elsewhere, with no residual raw-embedding constraint.
Finally,
\[
 I_w^{(j)}(\phi)=0
 \ \Longrightarrow\ \phi(t,\theta)=-wt+C,\qquad
 1=\int_0^\infty e^{(-wt+C)/2}\dd t=\frac2w e^{C/2}.
\]
Thus $C=2\log(w/2)$, proving the stated unique zero.
\end{proof}

\section{Limits of quantum cutting configurations}\label{sec:ends}\label{sec:stability}

The marginal LDPs for $A$ and $B$ are now available. To identify
the joint rate of $(A,B,M)$, we study convergent sequences of actual
cutting configurations whose limits have $I(A)+K(B)<\infty$.
This section passes the field change-of-coordinate identity and
boundary length matching to the limit. Section~\ref{sec:welding}
then proves coverage and uniqueness of the full configuration.

For the boundary extension lemma, use the following deterministic
setup, with $w>0$ denoting a component weight and $W$ the ambient weight.
The functions $\phi,\Phi$ are the corresponding real-valued field
profiles, locally $H^1$ up to bounded boundary intervals, and
\begin{equation}\label{eq:finitecutsetup}
 \begin{gathered}
 \phi(z)=-w\Rea z+u(z),\qquad
 \Phi(z)=-W\Rea z+U(z),\\
 \int_{\Sstrip}|\nabla u|^2\dd^2z<\infty,\qquad
 \int_{\Sstrip}|\nabla U|^2\dd^2z<\infty.
 \end{gathered}
\end{equation}
These energy bounds follow from finiteness of the field rates.
The candidate cutting map
$F:\Sstrip\to\Sstrip$ is holomorphic and injective and satisfies
\begin{equation}\label{eq:endscutidentity}
 \phi(z)=\Phi(F(z))+2\log|F'(z)|
 \quad\text{for a.e. }z\in\Sstrip.
\end{equation}
Its image is a component domain inside the ambient strip; we do
not assume that it fills the strip. In the application,
$(\phi,\Phi,F,w,W)$ stands for
$(\phi_i,\Phi,F_i,W_i,W)$.

\subsection{Boundary values and preservation of the ends}

The cutting identity is initially an identity inside the strip. To
pass boundary length matching to a limit, we must be able to evaluate
$F$ on the other boundary line, which will become the seam. Later,
when identifying the boundaries of the two compactified strips, we also need continuity at
their common endpoints. The following lemma supplies these boundary
values. Its endpoint conclusion says that approaching a strip end
from the interior or along either boundary gives the same ambient
endpoint as approaching along the known outer boundary.

\begin{lemma}[Continuous extension to the compactified strip]\label{lem:ends}
In the setup \eqref{eq:finitecutsetup}--\eqref{eq:endscutidentity},
fix an outer boundary index $j\in\{0,\pi\}$. Suppose $F$ extends
continuously to that boundary line, with
\[
 F(t+ij)=b(t)+ij,\qquad t\in\R,
\]
where $b:\R\to\R$ is an increasing homeomorphism. Then $F$
extends continuously to
\[
 F:\overline{\mathcal S}^{\,*}\longrightarrow
       \overline{\mathcal S}^{\,*},\qquad
 F(-\infty)=-\infty,\quad F(+\infty)=+\infty.
\]
Continuity is understood in the metric $d_\chi$ defined in the
introduction. Thus the conclusion supplies boundary values on
both strip boundary lines and preserves the two marked ends.
\end{lemma}
\begin{proof}
Choose a holomorphic logarithm of $F'$ and set
\[
 L=\log F'-WF/2+wz/2.
\]
The cutting identity cancels the two linear profiles:
\[
 2\Rea L=2\log|F'|-W\Rea F+w\Rea z=u-U\circ F.
\]
Thus $L$ has finite Dirichlet energy. This is the analytic control
from which both boundary continuity and the end estimate will follow.
Conformal invariance gives
\begin{equation}\label{eq:Lenergy}
 \int_{\Sstrip}|L'|^2
 \le\tfrac12\int_{\Sstrip}|\nabla u|^2
      +\tfrac12\int_{F(\Sstrip)}|\nabla U|^2<\infty.
\end{equation}
Schwarz reflection of the proper outer map makes $F$ analytic there
with positive derivative. Thus $L$ has constant imaginary part on
that line and reflects after subtraction of this constant.

\emph{Boundary extension.}
Put $G=e^{-W(F-ij)/2}$ and $q(t)=G(t+ij)>0$. On finite source
rectangles, $\log G'=L-wz/2+\text{constant}$ has finite Dirichlet
energy. Its exponential belongs to every finite $L^p$ by the
two-dimensional exponential inequality. Hence $G$ extends in
$W^{1,p}$, $p>2$, and continuously to bounded boundary intervals.
At zeros of the boundary extension of $G$, $\Rea F\to+\infty$;
since the imaginary part remains between $0$ and $\pi$, this is
convergence to the ambient end $+\infty$ in the compactified strip.
At nonzeros take a local logarithm to recover the continuous boundary
value of $F$. This proves continuity at each point of the two
boundary lines, with values in $\overline{\mathcal S}^{\,*}$.

\emph{Control along the known outer boundary.}
Properness on the outer line gives $q(+\infty)=0$, $q(-\infty)=\infty$
and $q'<0$. Reflection and mean-value bounds on fixed disks imply
$L'(t+ij)\to0$ as $t\to\pm\infty$, by \eqref{eq:Lenergy}. Since
\[
 \frac{d}{dt}\log(-q'(t))=-w/2+o(1),
\]
integration of the positive function $-q'$ at either end gives
$q'(t)/q(t)\to-w/2$. For example, at $+\infty$ compare
$-q'(t+s)/[-q'(t)]$ with $e^{-(w/2\pm\varepsilon)s}$ and integrate;
at $-\infty$ integrate backwards from a fixed time. Both comparisons
give the same limit.

\emph{Control across the whole strip width.}
On a fixed rectangle translated by $t$, Poincare's inequality and
reflection near the anchor $t+ij$ give
\[
 L(t+\cdot)-L(t+ij)\longrightarrow0\quad\text{in }H^1.
\]
To justify the point anchor, first subtract the rectangle mean; a
mean-value bound on a fixed reflected disk at the anchor controls its
value by the $H^1$ norm. The Dirichlet energy on enlarged translated
rectangles tends to zero. Exponential integrability now yields
$e^{L(t+\cdot)-L(t+ij)}\to1$ in every finite $L^p$.
Consequently
\[
 \frac{G'(t+\zeta)}{q(t)}\longrightarrow
       -\frac w2e^{-w(\zeta-ij)/2}\quad\text{in }L^p.
\]
Use $G(t+ij)/q(t)=1$ and $p>2$. The anchored $W^{1,p}$ inequality
and Sobolev embedding give uniform convergence on the closed rectangle:
\[
 \frac{G(t+\zeta)}{q(t)}\longrightarrow e^{-w(\zeta-ij)/2}.
\]
The right side never vanishes. The continuous logarithm anchored at
$\zeta=ij$ therefore converges to $-w(\zeta-ij)/2$, even when the
exponential winds through more than one full turn. Returning to $F$
gives the intermediate estimate, for each fixed $R<\infty$,
\begin{equation}\label{eq:endasympt}
 \sup_{|x|\le R,\,0\le y\le\pi}
 \left|F(t+x+iy)-F(t+ij)-\frac wW(x+i(y-j))\right|\to0
 \quad(t\to\pm\infty).
\end{equation}
In particular $F(t+iy)-F(t+ij)$ is uniformly bounded over
$0\le y\le\pi$ for large $|t|$. Since $F(t+ij)=b(t)+ij$
and $b(t)\to\pm\infty$ at the respective ends, we obtain
\[
 \sup_{0\le y\le\pi}d_\chi(F(t+iy),+\infty)\to0
 \quad(t\to+\infty),
\]
\[
 \sup_{0\le y\le\pi}d_\chi(F(t+iy),-\infty)\to0
 \quad(t\to-\infty).
\]
These are continuity at the two marked strip ends.
Together with the extension along both boundary lines already proved, they
give the continuous extension of $F$ to
$\overline{\mathcal S}^{\,*}$.
\end{proof}

\begin{remark}
The developing map $G$ is not assumed globally injective. Its possible
winding is why a planar convex-sector argument is unsuitable for large
weights. The anchored logarithm in the last step is essential.
\end{remark}

\subsection{Passing the cutting identities to the limit}

\begin{proposition}[Classical cutting and length matching in the limit]
\label{prop:closed}
Let $\kappa_n\downarrow0$ and consider actual quantum cutting
configurations in the embeddings of Subsection~\ref{subsec:ABM}.
Suppose $(A_n,B_n)\to(A,B)$ in $\mathcal E\times\mathcal E^2$,
with $I(A)+K(B)<\infty$, and $F_{i,n}\to F_i$ in $\mathcal F$.
Then the limits $F_i$ are nonconstant
injective holomorphic maps $\mathcal S\to\mathcal S$
with disjoint images, continuous in the two-end compactification,
and with increasing homeomorphisms on their prescribed outer lines.
They satisfy
\begin{equation}\label{eq:classicalcutlimit}
 \phi_i=\Phi\circ F_i+2\log|F_i'|\quad\text{a.e. on }\Sstrip
\end{equation}
and, writing $q_1=\ell_{1,\pi}^{-1}$ and
$q_2=\ell_{2,0}^{-1}$, the classical boundary length matching identity
\begin{equation}\label{eq:classicallengthmatching}
 \chi(F_1(q_1(l)+i\pi))=\chi(F_2(q_2(l)))
 \qquad(l>0).
\end{equation}
\end{proposition}
\begin{proof}
Write $\Phi_n,\phi_{i,n}$ for the fields and $\Phi,\phi_i$ for
their limits. Set $j_1=0,j_2=\pi$ for the outer boundary indices
and $k_1=\pi,k_2=0$ for the seam indices. The assumed enhanced-data
convergence includes distributional field convergence, strong local
$H^{-\sigma}$ convergence of measures, and local uniform convergence
of cumulative lengths. Compactness gives the assumed map subsequence,
but its limits could initially be constant. We first rule this out,
then pass the field identity and the boundary matching separately.

\emph{1. Sobolev convergence of the boundary maps.}
For this estimate only, use the bounded representatives
\[
 f_{i,n}=\chi\circ F_{i,n}\circ\chi^{-1},\qquad
 f_i=\overline\chi\circ F_i\circ\chi^{-1}.
\]
They converge locally uniformly on $\mathbb D$ by definition of
the map topology.
Injectivity and the area formula give
\[
 \int_\mathbb D|\nabla f_{i,n}|^2\dd^2z
 =2\int_\mathbb D|f_{i,n}'|^2\dd^2z
 =2\operatorname{Area}(f_{i,n}(\mathbb D))
 \le2\operatorname{Area}(\mathbb D).
\]
Since $\mathbb D$ is bounded, the maps are bounded in $H^1(\mathbb D)$.
Local uniform convergence identifies their weak $H^1$ limit as $f_i$.
The trace operator is bounded into $H^{1/2}(\partial\mathbb D)$,
and the embedding into $H^\sigma(\partial\mathbb D)$ is compact for
$\sigma<1/2$. Therefore
\begin{equation}\label{eq:strongmaptrace}
 \operatorname{Tr}f_{i,n}\longrightarrow\operatorname{Tr}f_i
 \quad\text{strongly in }H^\sigma(\partial\mathbb D).
\end{equation}
Uniqueness of the limit makes this convergence valid for the full
sequence. On a bounded source boundary interval, composition with
the smooth chart $t\mapsto\chi(t+ij)$ and multiplication by a
smooth cutoff preserve this convergence.

\emph{2. Outer maps and nondegeneracy.}
The purpose of this step is to rule out constant limits. A locally
uniform limit of injective holomorphic maps is either injective
holomorphic or constant. A constant limit would lose the component
parametrization, and a constant on $\partial\mathbb D$ would not even
allow us to return to strip coordinates using $\chi^{-1}$.
The boundary length change-of-coordinate formula gives
\[
 \ell_{i,j_i,n}(t)
 =\ell_{\mathrm{amb},j_i,n}(\Rea F_{i,n}(t+ij_i)),
\]
Define the candidate limiting outer map by
\begin{equation}\label{eq:outerlimitmap}
 b_i(t)=\ell_{\mathrm{amb},j_i}^{-1}(\ell_{i,j_i}(t)).
\end{equation}
At this point $b_i$ is defined solely by the limiting lengths;
we have not yet identified it with a boundary value of $F_i$.
The same change-of-coordinate formula gives
\[
 \Rea F_{i,n}(t+ij_i)
 =\ell_{\mathrm{amb},j_i,n}^{-1}(\ell_{i,j_i,n}(t))
 \longrightarrow b_i(t).
\]
This convergence is uniform on compact $t$-intervals: the limiting
cumulative functions are continuous strictly decreasing bijections,
and their inverses converge uniformly on compact length intervals
by the same bracketing argument as \eqref{eq:lengthrootbracket}.
The composition of the two decreasing bijections in
\eqref{eq:outerlimitmap} is increasing and satisfies
\[
 b_i(-\infty)=-\infty,\qquad b_i(+\infty)=+\infty.
\]
To connect this boundary convergence to the interior limit, fix a
bounded interval $J$. With $\gamma_i(t)=\chi(t+ij_i)$, we have
two limits for the same boundary functions:
\[
 \begin{aligned}
 f_{i,n}(\gamma_i(t))
   &\longrightarrow\chi(b_i(t)+ij_i)
       &&\text{uniformly for }t\in\overline J,\\
 (\operatorname{Tr}f_{i,n})(\gamma_i(t))
   &\longrightarrow(\operatorname{Tr}f_i)(\gamma_i(t))
       &&\text{in }L^2(J),
 \end{aligned}
\]
where the second convergence follows from \eqref{eq:strongmaptrace}.
Thus
\[
 (\operatorname{Tr}f_i)(\gamma_i(t))
       =\chi(b_i(t)+ij_i)\quad\text{for a.e. }t\in J.
\]
If $f_i\equiv c$, its Sobolev trace would also be $c$ a.e.
But $t\mapsto\chi(b_i(t)+ij_i)$ is continuous and injective,
so it cannot equal a constant a.e. on an interval. This contradiction
rules out every constant limit, including constants on the unit circle.
The theorem on locally uniform limits of injective holomorphic maps
then implies that $f_i$ is injective. The open mapping theorem puts
its image inside $\mathbb D$, giving $F_i$ as stated. The images of
$F_1,F_2$ are disjoint: a point in both limiting images would, by
local inverse convergence (or Rouch\'e's theorem), belong to both
approximating images for all large $n$, which is impossible.
The identified continuous outer traces give continuous outer boundary
values; Schwarz reflection across each bounded interval of the outer boundary gives
the analytic extension there.
Consequently these boundary values satisfy
\[
 F_i(t+ij_i)=b_i(t)+ij_i,\qquad t\in\mathbb R.
\]

\emph{3. The interior identity and boundary extensions.}
Step 2 establishes a conformal map, but does not identify its relation
to the limiting fields. We must still prove
\eqref{eq:classicalcutlimit}. This is also the missing hypothesis
needed to apply Lemma~\ref{lem:ends}: convergence inside the strip
alone does not give continuous values on the seam or at the ends.
For $\zeta\in C_c^\infty(\Sstrip)$, change of variables gives
\[
 \langle\Phi_n\circ F_{i,n},\zeta\rangle
 =\langle\Phi_n,T_{i,n}\zeta\rangle,\qquad
 T_{i,n}\zeta(w)=\zeta(F_{i,n}^{-1}(w))
                         |(F_{i,n}^{-1})'(w)|^2,
\]
with the latter test extended by zero outside $F_{i,n}(\Sstrip)$.
The support of $T_{i,n}\zeta$ lies in
$F_{i,n}(\operatorname{supp}\zeta)$. Local uniform convergence to
the nondegenerate $F_i$ puts these supports in a common compact
subset of the ambient strip. Convergence of the holomorphic maps
and their local inverses, with all derivatives on compact sets,
then gives $T_{i,n}\zeta\to T_i\zeta$ in $C_c^\infty$.
Distributional convergence of $\Phi_n$
then gives
\[
 \langle\Phi_n,T_{i,n}\zeta\rangle\to
 \langle\Phi,T_i\zeta\rangle,
 \qquad
 (2+\kappa_n/2)\log|F_{i,n}'|\to2\log|F_i'|
 \quad\text{locally smoothly}.
\]
Passing to the limit in \eqref{eq:quantumcut} proves
\eqref{eq:classicalcutlimit}. Lemma~\ref{lem:ends} now supplies
continuous boundary extensions of $\chi\circ F_i$, including at
both ends. These continuous extensions represent their Sobolev traces.

\emph{4. Boundary errors integrated against the varying measures.}
Interior convergence has not yet passed the seam identity to the
limit. We first obtain an integrated boundary error estimate in source
coordinates, and then change to length coordinates in the next step.
On the seam write
\[
 B_{i,n}(t)=\chi(F_{i,n}(t+ik_i)),\qquad
 B_i(t)=\chi(F_i(t+ik_i)),\qquad
 m_{i,n}=\nu_{i,k_i,n},\quad m_i=\nu_{i,k_i}.
\]
Fix a finite interval $K$ and smooth nested nonnegative cutoffs
$\zeta,\zeta_1,\zeta_2$, with $\zeta=1$ on $K$ and the larger
cutoffs equal to one on all supports needed below. The assumed
measure convergence gives
\[
 \sup_n\|\zeta_1m_{i,n}\|_{H^{-\sigma}}<\infty,
\]
whereas \eqref{eq:strongmaptrace} gives
$\|\zeta_2(B_{i,n}-B_i)\|_{H^\sigma}\to0$. Since absolute value
is bounded on fractional $H^\sigma$, Sobolev duality yields
\begin{equation}\label{eq:pairing}
 \begin{split}
 \int_K|B_{i,n}-B_i|\dd m_{i,n}
 &\le\int\zeta|B_{i,n}-B_i|\dd m_{i,n}\\
 &\le C\|\zeta_1m_{i,n}\|_{H^{-\sigma}}
          \|\zeta_2(B_{i,n}-B_i)\|_{H^\sigma}
 \longrightarrow0.
 \end{split}
\end{equation}
The pairing is legitimate for these positive measures: a positive
measure in $H^{-\sigma}$ does not charge sets of zero
$H^\sigma$ capacity and integrates the quasi-continuous Sobolev
representative. Here the boundary maps are continuous, so those
representatives are exactly the boundary values being used.

\emph{5. Change to length coordinates and pass the matching identity.}
Fix $0<a<b<\infty$ and put
$q_{i,n}=\ell_{i,k_i,n}^{-1}$. Local uniform cumulative convergence
and strict monotonicity give
\begin{equation}\label{eq:inversecompactlength}
 \sup_{l\in[a,b]}|q_{i,n}(l)-q_i(l)|\longrightarrow0,
 \qquad q_{i,n}([a,b])\cup q_i([a,b])\subset K_i
\end{equation}
for a fixed compact source interval $K_i$ and all large $n$.
The decreasing cumulative coordinate sends its own boundary measure
to Lebesgue length. Explicitly, for $0<c<d$,
\[
 m_{i,n}([q_{i,n}(d),q_{i,n}(c)])=d-c,
\]
and hence, for every nonnegative Borel $H$,
\begin{equation}\label{eq:lengthsubstitution}
 \int_a^b H(q_{i,n}(l))\dd l
 =\int_{q_{i,n}(b)}^{q_{i,n}(a)}H(t)\dd m_{i,n}(t).
\end{equation}
Atomlessness makes endpoint conventions irrelevant. No derivative
of $q_{i,n}$ or density of $m_{i,n}$ is needed.

Define the seam maps parametrized by length by
$C_{i,n}(l)=B_{i,n}(q_{i,n}(l))$ and $C_i(l)=B_i(q_i(l))$.
Split their error into two terms:
\[
 \begin{split}
 \int_a^b|C_{i,n}-C_i|\dd l
 &\le\int_a^b|(B_{i,n}-B_i)(q_{i,n}(l))|\dd l\\
 &\quad+\int_a^b|B_i(q_{i,n}(l))-B_i(q_i(l))|\dd l\\
 &\le\int_{K_i}|B_{i,n}-B_i|\dd m_{i,n}
       +(b-a)\,\omega_{i,K_i}\!\left(
          \sup_{[a,b]}|q_{i,n}-q_i|\right)
 \longrightarrow0.
 \end{split}
\]
Here $\omega_{i,K_i}$ is the modulus of continuity of $B_i$ on
$K_i$. The first term tends to zero by \eqref{eq:pairing}, and
the second by \eqref{eq:inversecompactlength} and continuity of
$B_i$. Exact quantum matching says $C_{1,n}(l)=C_{2,n}(l)$.
Therefore
\[
 \int_a^b|C_1(l)-C_2(l)|\dd l
 \le\|C_1-C_{1,n}\|_{L^1([a,b])}
    +\|C_{2,n}-C_2\|_{L^1([a,b])}\longrightarrow0.
\]
Thus $C_1=C_2$ almost everywhere on $[a,b]$. Both limiting maps
are continuous, so equality holds everywhere there. Since
$[a,b]\subset(0,\infty)$ was arbitrary, this proves
\eqref{eq:classicallengthmatching} for every $l>0$.
\end{proof}
\section{Classical conformal welding and uniqueness}\label{sec:welding}

We prove that fixed finite-rate component data $B$ determine at
most one limiting pair $(A,M)$ with $I(A)<\infty$. Section~\ref{sec:ends}
already gives the cutting and boundary length identities. We first
show that the limiting cutting maps fill the ambient strip, meet
along a Jordan crosscut, and determine the normalized ambient field
uniquely. We then identify this crosscut with the separately recorded
Hausdorff limit of the curves and identify the geometric maps in $M$.

The geometric argument joins the two source strips by their length
correspondence and constructs conformal charts across the join,
following the classical method of \cite[Section 3.1]{RS}.
\begin{definition}[Bounded and vanishing mean oscillation]
\label{def:VMO}
For a real locally integrable function $a$ and a bounded interval
$J$, write
\[
 a_J=\frac1{|J|}\int_J a(t)\dd t,\qquad
 \operatorname{osc}_J(a)=\frac1{|J|}\int_J|a-a_J|\dd t.
\]
On a bounded interval $J_0$, the BMO seminorm is
\[
 \|a\|_{\mathrm{BMO}(J_0)}
 =\sup_{J\subset J_0}\operatorname{osc}_J(a).
\]
The function has locally bounded mean oscillation if this seminorm
is finite on every compact interval. It has locally vanishing mean
oscillation, written $a\in\mathrm{VMO}_{\mathrm{loc}}(\R)$, if
\[
 \lim_{r\downarrow0}
 \sup_{\substack{J\subset K\\0<|J|\le r}}
 \operatorname{osc}_J(a)=0
 \quad\text{for every compact interval }K\subset\R.
\]
In particular, around any point a VMO function has arbitrarily
small BMO seminorm on a sufficiently small interval.
\end{definition}

\begin{lemma}[John--Nirenberg inequality]\label{lem:JN}
There are universal constants $c,C>0$ such that, for any bounded
interval $J_0$ and real function $a$ with
$\|a\|_{\mathrm{BMO}(J_0)}\le\varepsilon$, where $\varepsilon>0$,
\[
 \frac{|\{t\in J:|a(t)-a_J|>\lambda\}|}{|J|}
 \le C e^{-c\lambda/\varepsilon},
 \qquad J\subset J_0,\quad\lambda>0.
\]
Consequently, for $0<p<c/\varepsilon$,
\[
 \frac1{|J|}\int_J e^{p|a-a_J|}\dd t
 \le 1+\frac{Cp}{c/\varepsilon-p},\qquad J\subset J_0.
\]
In particular, if $\varepsilon\le c/4$, then
$|J|^{-1}\int_J e^{2|a-a_J|}\dd t\le 1+C$ uniformly in $J$.
\end{lemma}
\begin{proof}
The first assertion is the classical John--Nirenberg inequality \cite{JohnNirenberg1961};
see \cite[Section 2.1]{VW}. Its exponential-moment consequence
follows by integrating the tail:
\[
 \frac1{|J|}\int_J e^{p|a-a_J|}\dd t
 =1+p\int_0^\infty e^{p\lambda}
   \frac{|\{t\in J:|a(t)-a_J|>\lambda\}|}{|J|}\dd\lambda
 \le 1+Cp\int_0^\infty e^{-(c/\varepsilon-p)\lambda}\dd\lambda.
\]
\end{proof}

We call
$(\Phi,F_1,F_2)$ a \emph{candidate welding} of the component fields
$(\phi_1,\phi_2)$ if all three fields have finite wedge energy,
the maps $F_i:\Sstrip\to\Sstrip$ are injective holomorphic maps
with disjoint images, and
\[
 \phi_i=\Phi\circ F_i+2\log|F_i'|,\qquad
 F_1(\ell_{1,\pi}^{-1}(l)+i\pi)
 =F_2(\ell_{2,0}^{-1}(l))\quad(l>0).
\]
Boundary values are interpreted in the two-end compactification.
On the outer lines, require $F_1(t)\in\R$ and
$F_2(t+i\pi)\in\R+i\pi$, with each boundary map an increasing
homeomorphism onto that line. These are precisely the properties
proved for a limit in Proposition~\ref{prop:closed}. A candidate is
\emph{normalized} when the component and ambient fields use the
boundary-distance embeddings fixed above; in particular $F_1(0)=0$.

\begin{theorem}[Coverage and uniqueness of the welding]\label{thm:glue}
Every candidate welding satisfies
\[
\Sstrip=F_1(\Sstrip)\,\dot\cup\,(\eta_*\cap\Sstrip)
                       \,\dot\cup\,F_2(\Sstrip),
\]
where $\eta_*$ is a Jordan crosscut joining the two ends, locally a
quasiarc and of zero planar area. For fixed component fields, there
is at most one normalized candidate welding.
\end{theorem}
\begin{proof}
Identify the seam boundary arcs of the two compactified source strips using the
classical length correspondence. Lemma~\ref{lem:classcdf} makes it a
homeomorphism of closed arcs including the endpoints. The resulting
topological space $X$ is a closed disk; its boundary consists of the
two outer arcs.

We next show why the length correspondence is locally
quasisymmetric; this supplies the conformal charts across the seam.
Write the seam length densities as
\[
 \rho_1(t)=e^{a_1(t)},\quad a_1(t)=\tfrac12\phi_1(t,\pi),
 \qquad
 \rho_2(t)=e^{a_2(t)},\quad a_2(t)=\tfrac12\phi_2(t,0).
\]
Finite wedge energy and the Sobolev trace theorem give
$a_i\in H^{1/2}_{\mathrm{loc}}(\R)$. The affine part of each field
is smooth and causes no difficulty on bounded intervals.

We verify the VMO property of Definition~\ref{def:VMO} for each
$a=a_i$. Cauchy--Schwarz gives
\[
 \left(\frac1{|J|}\int_J|a-a_J|\dd t\right)^2
 \le \frac1{|J|^2}\int_J\int_J|a(t)-a(s)|^2\dd t\dd s
 \le \int_J\int_J
       \frac{|a(t)-a(s)|^2}{|t-s|^2}\dd t\dd s.
\]
On a fixed bounded interval the last integrand is integrable by
the $H^{1/2}$ assumption. Absolute continuity of its integral
therefore shows that
\[
 \sup_{\substack{J\subset K\\0<|J|\le r}}
 \frac1{|J|}\int_J|a-a_J|\dd t\longrightarrow0
 \quad(r\downarrow0)
\]
for every compact interval $K$. Near any fixed seam point we may
therefore choose $J_0$ with
$\|a\|_{\mathrm{BMO}(J_0)}\le c/4$, where $c$ is the constant
in Lemma~\ref{lem:JN}. That lemma, with $p=2$, gives a uniform
constant $C_0$ such that
\[
 \frac1{|J|}\int_J e^{2|a-a_J|}\dd t\le C_0,
 \qquad J\subset J_0.
\]
Thus $\rho=e^a$ and its reciprocal obey the local $A_2$ bound
\[
 \left(\frac1{|J|}\int_J\rho\dd t\right)
 \left(\frac1{|J|}\int_J\rho^{-1}\dd t\right)
 \le C_0^2.
\]
For clarity, an $A_\infty$ estimate also follows directly from
Cauchy--Schwarz and Jensen: for measurable $E\subset J$,
\[
 \frac{\int_E\rho\dd t}{\int_J\rho\dd t}
 \le C_0^{1/2}\left(\frac{|E|}{|J|}\right)^{1/2}.
\]
The $A_2$ bound implies comparability of the masses of adjacent
intervals of equal length. Indeed, for $I\subset J$,
\[
 \frac{\int_I\rho\dd t}{\int_J\rho\dd t}
 \ge C_0^{-2}\left(\frac{|I|}{|J|}\right)^2,
\]
by Cauchy--Schwarz applied to $1=\rho^{1/2}\rho^{-1/2}$ on $I$.
Apply this to the two halves of $J$. Consequently the increasing
primitive $L(t)=\int_{t_0}^t\rho(s)\dd s$ satisfies
\[
 Q^{-1}\le
 \frac{L(t+r)-L(t)}{L(t)-L(t-r)}\le Q
\]
whenever the intervals lie in $J_0$, for a fixed $Q\ge1$.
This is local quasisymmetry.

Applying this argument to both densities, their increasing primitives
$L_1,L_2$ are locally quasisymmetric. The matching map
$h=\ell_{2,0}^{-1}\circ\ell_{1,\pi}$ satisfies
$L_2(h(t))=L_1(t)+c$ locally, for a constant $c$. Since inverses
and compositions of quasisymmetric maps are quasisymmetric,
$h$ is locally quasisymmetric as well. This is the local form of
the isometric welding argument in \cite[Section 2.1 and
Proposition 2.5]{VW}. Local quasisymmetric welding
provides conformal charts across the seam in which it is a quasiarc.
These charts and the original source charts define a Riemann surface
structure on $X^\circ$.

The bounded maps $\chi\circ F_i$ agree on the seam, so assemble into a
continuous map $f:X\to\overline{\mathbb D}$. Continuity at the endpoints
follows from Lemma~\ref{lem:ends}. Each restriction has Dirichlet
energy at most $2\operatorname{Area}(\mathbb D)$. In a local welding chart,
straighten the quasiarc by a quasiconformal map. Dirichlet energies
remain finite on both sides. Equality of the continuous traces and
the Sobolev gluing lemma across a straight line give $W^{1,2}$ across
the seam; changing back preserves this property. In the conformal
chart $\bar\partial f=0$ a.e., since the seam has area zero. Weyl's
lemma proves holomorphicity across the seam. The map is nonconstant,
so the open mapping theorem puts all of $X^\circ$ in $\mathbb D$.

On $\partial X$, the two proper outer maps combine into an
orientation-preserving homeomorphism onto $\partial\mathbb D$. The map
is proper over interior compact sets and has degree one. For any
interior value, its finitely many holomorphic preimages have positive
integer multiplicities whose sum is one. Thus $f$ is a conformal
bijection $X^\circ\to\mathbb D$, with a homeomorphism of compactifications.
The image of the identified boundary arc is consequently a Jordan crosscut. Its
interior is locally a quasiarc and has area zero. This proves full
coverage without any finite-total-area input.

For uniqueness, fix the component fields. Their cumulative lengths
fix the identification $h=\ell_{2,0}^{-1}\circ\ell_{1,\pi}$,
and hence the welded Riemann surface $X^\circ$ with the conformal
structure just constructed. If $F_i$ and $\widetilde F_i$ are
two configurations, they induce conformal bijections
$T,\widetilde T:X^\circ\to\Sstrip$. Consequently
\[
 H=\widetilde T\circ T^{-1}:\Sstrip\to\Sstrip,\qquad
 \widetilde F_i=H\circ F_i.
\]
The automorphism $H$ preserves the two ends individually, so
$H(z)=z+c$ for some real constant $c$. In the prescribed embeddings,
the first component and ambient lower boundary length-one marks both have
coordinate zero. Therefore
\[
 F_1(0)=\widetilde F_1(0)=0\quad\Longrightarrow\quad c=0.
\]
The maps and their common seam are thus identical. The cutting
identity determines the ambient field on each image by
\[
 \Phi(F_i(z))=\phi_i(z)-2\log|F_i'(z)|.
\]
Since the images cover the strip apart from a zero-area seam,
this determines the ambient field almost everywhere. Its classical
measures and cumulative lengths are then fixed as well. This proves
uniqueness of the normalized ambient data and cutting maps.
\end{proof}
The preceding theorem identifies a classical crosscut from the
limiting cutting maps. We must still show that this crosscut agrees
with the curve coordinate recorded in $M$, whose convergence is
measured in Hausdorff distance.

\begin{lemma}[Identification of the limiting curve and geometric maps]
\label{lem:trace}
Let $\kappa_n\downarrow0$ and let
\[
(A_n,B_n,M_n)\longrightarrow(A,B,M)
\]
be a sequence of actual quantum cutting configurations, with
$I(A)+K(B)<\infty$. Let $\eta_*$ be the classical crosscut determined
by the limiting maps $F_1,F_2$ in Theorem~\ref{thm:glue}, including
its two marked ends. Then the curve coordinate $\eta$ of $M$
equals $\eta_*$. Moreover, the geometric maps $G_{i,n}$ derived
from the actual configurations converge locally uniformly to
the geometrically normalized maps $G_i$ onto $F_i(\mathcal S)$.
The recovery formula \eqref{eq:derivedG} is continuous at each
finite-rate welding configuration.

Consequently, for fixed finite-rate component data $B$, there is
at most one pair $(A,M)$ obtainable as such a limit with
$I(A)<\infty$.
\end{lemma}

\begin{proof}
Proposition~\ref{prop:closed} and Theorem~\ref{thm:glue} identify
the limiting domains
\[
D_i=F_i(\mathcal S),
\qquad
\mathcal S
=D_1\,\dot\cup\,(\eta_*\cap\mathcal S)\,\dot\cup\,D_2.
\]
It remains to identify the separately recorded Hausdorff limit
of the curves.

For this argument, use disk coordinates and write
\[
f_{i,n}=\chi\circ F_{i,n}\circ\chi^{-1},
\qquad
f_i=\chi\circ F_i\circ\chi^{-1},
\]
\[
E_{i,n}=f_{i,n}(\mathbb D),\qquad
E_i=f_i(\mathbb D),\qquad
\widehat\eta=\overline\chi(\eta),\qquad
\Gamma=\overline\chi(\eta_*).
\]
Thus $\overline\chi(\eta_n)\to\widehat\eta$ in Hausdorff distance,
and our goal is to prove $\widehat\eta=\Gamma$.

\emph{No extra points inside the disk.}
Local uniform convergence to the injective maps $f_i$ gives
\[
C\Subset E_i
\quad\Longrightarrow\quad
C\subset E_{i,n}
\quad\text{for all sufficiently large }n.
\]
Indeed, this follows from Rouch\'e's theorem applied near the
preimages under $f_i$, followed by a finite covering argument.
A small closed disk contained in $E_i$ therefore eventually
avoids the approximating curves. Since
\[
\mathbb D=E_1\,\dot\cup\,(\Gamma\cap\mathbb D)\,\dot\cup\,E_2,
\]
we obtain
\[
\widehat\eta\cap\mathbb D\subset\Gamma.
\]

\emph{No extra points on the outer boundary.}
Fix $p\in\partial\mathbb D\setminus\{-1,1\}$. The outer boundary
homeomorphisms give a unique preimage $q$ of $p$ on the
corresponding source outer arc. Choose an analytic coordinate
$\psi$ mapping an upper half-disk
\[
H_r=\{z:|z|<r,\ \operatorname{Im}z>0\}
\]
into the source disk, with $\psi(0)=q$ and its diameter on that
outer arc. Set
\[
g_n=f_{i,n}\circ\psi,\qquad g=f_i\circ\psi.
\]
The local uniform convergence of the outer boundary maps proved
in Proposition~\ref{prop:closed} allows us to choose $r$ so that
\[
\sup_{x\in[-r,r]}|g_n(x)-p|\le\tfrac14
\]
for all sufficiently large $n$.

Let $\omega$ be the harmonic measure of the semicircular part of
$\partial H_r$. Since $|g_n-p|\le2$ in $H_r$, the subharmonic
maximum principle gives
\[
|g_n(z)-p|\le\tfrac14+2\omega(z).
\]
As $\omega(z)\to0$ when $z\to0$, there is a fixed $\rho>0$ such
that
\[
|g_n(z)-p|\le\tfrac12,\qquad |g_n(z)|\ge\tfrac12
\quad(z\in H_\rho)
\]
for all large $n$. Reflection across the target circle therefore
gives holomorphic extensions
\[
\widetilde g_n(z)=
\begin{cases}
g_n(z),&\operatorname{Im}z\ge0,\\[2mm]
1/\overline{g_n(\overline z)},&\operatorname{Im}z<0
\end{cases}
\]
to $|z|<\rho$, uniformly bounded by $2$.

These extensions are injective: the upper half maps inside the
unit disk, the lower half outside its closure, and the diameter
maps injectively to the circle. Their locally uniform limit is
the nonconstant holomorphic extension $\widetilde g$ of $g$,
and is therefore injective. Since $\widetilde g(0)=p$,
Rouch\'e's theorem on a smaller circle centered at zero gives
a fixed neighborhood $V$ of $p$ contained in
$\widetilde g_n(\{|z|<\rho\})$ for all large $n$.
A point of $V\cap\mathbb D$ has its preimage in the upper half,
so
\[
V\cap\mathbb D\subset E_{i,n}.
\]
After shrinking $V$ to avoid the marked endpoints, it therefore
avoids the entire approximating curve. Thus $p\notin\widehat\eta$.
Together with the preceding step, this proves $\widehat\eta\subset\Gamma$.

\emph{Every point of the classical crosscut is attained.}
Suppose $p\in\Gamma\cap\mathbb D$ but $p\notin\widehat\eta$. A sufficiently
small disk around $p$ then avoids the approximating curves for
all large $n$. Since $p$ lies on the common boundary of $E_1$
and $E_2$, this disk contains points $x_i\in E_i$.
Interior compact containment gives $x_i\in E_{i,n}$ for all
large $n$. This is impossible: a connected disk avoiding the
approximating crosscut cannot contain points from both of its
complementary components. Finally, the marked endpoints belong
to every approximating curve. Hence $\Gamma\subset\widehat\eta$, proving
\[
\widehat\eta=\Gamma,\qquad \eta=\eta_*.
\]

\emph{Recovery of the geometric maps.}
Use the indices and formula in \eqref{eq:derivedG}:
\[
 s_{i,n}=\ell_{i,j_i,n}^{-1}
             \bigl(\ell_{\mathrm{amb},j_i,n}(0)\bigr),\qquad
 s_i=\ell_{i,j_i}^{-1}\bigl(\ell_{\mathrm{amb},j_i}(0)\bigr).
\]
Local uniform convergence of the cumulative lengths, positivity
of the limiting length at zero, and strict monotonicity give
$s_{i,n}\to s_i$ by the bracketing argument
\eqref{eq:lengthrootbracket}. Therefore
\[
G_{i,n}(z)=F_{i,n}(z+s_{i,n})
\longrightarrow F_i(z+s_i).
\]
The latter is precisely the geometrically normalized map onto
$D_i$. This also proves continuity of the recovery formula at
finite-rate data: it uses only convergence of cumulative lengths
and interior convergence of the cutting maps with bounded source
translations. Off the actual data, extend the inverse selection
measurably as in the proof of Theorem~\ref{thm:wedge}; that selection
remains continuous at the strictly monotone finite-rate lengths.

Theorem~\ref{thm:glue} determines the normalized ambient field
and cutting maps uniquely from $B$. The identification above
determines the recorded curve as well.
This proves uniqueness of the full limiting pair $(A,M)$.
\end{proof}
\section{Large deviations for the welded curve}\label{sec:curve}

We now use the marginal LDPs for $A$ and $B$ and the uniqueness of
the limiting pair $(A,M)$ proved in Section~\ref{sec:welding} to
identify the joint rate of $(A,B,M)$. Exponential tightness gives
subsequential joint LDPs; marginal contraction and uniqueness force
every such rate to equal $K(B)$ at the unique finite-rate welding
configuration. We then project to the interface and use independence
of the ambient wedge and SLE to identify its rate at the zero-energy
ambient field.

\subsection{The joint LDP for $(A,B,M)$}

\begin{proposition}[Joint cutting LDP]\label{prop:joint}
The tuple $(A_\kappa,B_\kappa,M_\kappa)$ satisfies a good LDP at
speed $\kappa^{-1}$ in
$\mathcal E\times\mathcal E^2\times\mathcal M$, with rate
$L$ given by \eqref{eq:ABMtarget}. For every $B$ with $K(B)<\infty$,
the normalized finite-rate welding configuration in that formula
exists and is unique.
\end{proposition}
\begin{proof}
Theorem~\ref{thm:wedge} gives good exponentially tight marginal
LDPs for $A_\kappa$ and $B_\kappa$, with rates $I$ and $K$.
The component rate is the sum of the two enhanced wedge rates
because the component surfaces are independent. The two map
coordinates lie in the compact space $\mathcal F^2$, and the curve
coordinate lies in the compact Hausdorff space
$\mathcal K_\chi(\overline{\mathcal S}^{\,*})$. Thus $\mathcal M$
is compact. For compact marginal sets $C_A,C_B$, the union bound
gives
\[
 \Pp\{(A_\kappa,B_\kappa,M_\kappa)
             \notin C_A\times C_B\times\mathcal M\}
 \le\Pp\{A_\kappa\notin C_A\}
       +\Pp\{B_\kappa\notin C_B\}.
\]
Consequently the joint laws are exponentially tight. The
subsequence theorem for exponentially tight laws on Polish spaces
\cite[Section 4.1]{DZ} implies that every sequence
$\kappa_n\downarrow0$ has a further subsequence with a good joint
LDP. Fix such a subsequence and denote its initially unknown rate
by $L$. Marginal contraction yields
\[
 \begin{gathered}
 I(A)=\inf_{B,M}L(A,B,M),\qquad
 K(B)=\inf_{A,M}L(A,B,M),\\
 L(A,B,M)\ge\max\{I(A),K(B)\}.
 \end{gathered}
\]

If $L(z)<\infty$ for $z=(A,B,M)$, every open neighborhood $U$ of
$z$ has positive probability at all sufficiently small parameters
along this subsequence: the LDP lower bound gives
\[
 \liminf_n\kappa_n\log\Pp\{(A_{\kappa_n},B_{\kappa_n},
                 M_{\kappa_n})\in U\}
 \ge-\inf_U L\ge-L(z)>-\infty.
\]
Choosing outcomes in shrinking neighborhoods therefore supplies
a sequence of actual configurations converging to $z$. This is
a deterministic selection of outcomes, not a claim that typical
random configurations converge to every finite-rate point.
Moreover $I(A)+K(B)<\infty$. Proposition~\ref{prop:closed},
Theorem~\ref{thm:glue}, and Lemma~\ref{lem:trace} show that $z$
is the normalized classical welding configuration determined by
$B$, including all three coordinates of $M$. In particular,
\[
 \#\{(A,M):L(A,B,M)<\infty\}\le1.
\]

If $K(B)<\infty$, choose approximate minimizers with
$L(A_n,B,M_n)\le K(B)+1/n$. They lie in the compact sublevel
$\{L\le K(B)+1\}$. A convergent subsequence and lower
semicontinuity give an attained minimum
\[
 L(A_*,B,M_*)=K(B).
\]
It is the unique finite-rate point above $B$. If $K(B)=\infty$,
every point above $B$ has infinite joint rate. This identifies
$L$ as \eqref{eq:ABMtarget} and proves existence of the configuration
for each finite-rate $B$.

Lemma~\ref{lem:trace} gives uniqueness among limits of actual
configurations independently of the subsequence. Hence every
subsequential rate has the same finite domain and values. The
subsequence characterization of the LDP gives the full good joint
LDP.
\end{proof}

\subsection{Projection and identification of the curve rate}
We have established that $(A_\kappa, B_\kappa, M_\kappa)$ satisfies
an LDP with rate function $L(A,B,M)$. Moreover, whenever
$L(A,B,M)<\infty$, we have $L(A,B,M)=K(B)$.
Projection therefore gives an LDP for the interface, with rate
\[
J(\eta)=\inf_{A,B,F_1,F_2}
L\bigl(A,B,(F_1,F_2,\eta)\bigr).
\]
It remains to identify this infimum with the cut-field energy
\[
E(\eta):=
I_{W_1}^{(0)}(\phi_1^\eta)
+
I_{W_2}^{(\pi)}(\phi_2^\eta).
\]
Independence of the ambient field and interface allows us to fix
the ambient datum at its zero-energy value. The cutting identities
then give $J(\eta)=E(\eta)$ whenever $J(\eta)<\infty$.
To obtain equality everywhere, we also prove that
$E(\eta)<\infty$ implies $J(\eta)<\infty$. For this step, the
boundary regularity lemma below shows that the deterministic
cut fields satisfy boundary length matching, so uniqueness
identifies their welding interface with $\eta$.

\begin{lemma}[Regularity along a boundary interval]
\label{lem:localregularity}
In the setup \eqref{eq:finitecutsetup}--\eqref{eq:endscutidentity},
let $I\subset\R$ be a bounded open interval and $k\in\{0,\pi\}$.
Suppose $F$ extends continuously to $I+ik$, the boundary map
$t\mapsto F(t+ik)$ is injective, and every value $F(t+ik)$ is
a point of $\mathbb C$, rather than one of the two compactified
strip ends. Thus the image of $I+ik$ is a curve parametrized
continuously and injectively by $t\in I$. Then:
\begin{enumerate}
 \item The boundary map $t\mapsto F(t+ik)$ is locally absolutely
 continuous, with nonzero derivative almost everywhere. Its image
 is locally rectifiable, and for each compact interval $J\Subset I$
 and every nonnegative Borel function $h$ on the image,
 \begin{equation}\label{eq:arclengthchange}
  \int_{F(J+ik)}h(\zeta)\dd s(\zeta)
  =\int_Jh(F(t+ik))|F'(t+ik)|\dd t.
 \end{equation}
 Here $\dd s$ is arclength and the boundary value of $F'$ is its
 almost everywhere trace.
 \item For every $p\in I+ik$, there are open sets
 $N_p,N'_p\subset\mathbb C$, containing $p$ and $F(p)$
 respectively, and a quasiconformal homeomorphism
 \[
 \widetilde F_p:N_p\longrightarrow N'_p,
 \qquad \widetilde F_p=F\ \text{on }N_p\cap\Sstrip.
 \]
 Thus $F$ extends across the boundary near $p$. The inverse
 $\widetilde F_p^{-1}$ maps the image boundary curve back to a
 straight interval and preserves local $H^1$ regularity, with
 a bound on Dirichlet energy depending on its distortion.
\end{enumerate}
\end{lemma}

A quasiconformal extension here means a planar homeomorphism
$\widetilde F\in W^{1,2}_{\rm loc}$ agreeing with $F$ on the
source side and satisfying
\[
 |\bar\partial\widetilde F|\le k_0|\partial\widetilde F|
 \quad\text{a.e.},\qquad k_0<1.
\]
Its distortion bound may depend on the neighborhood. The first
conclusion of the lemma concerns boundary lengths; the second
concerns gluing Sobolev functions across the image arc. Neither
conclusion is a convergence statement about approximating maps.
\begin{proof}
On a sufficiently small source half-disk $V$ the image of $F$ is
bounded. Equation \eqref{eq:Lenergy} and
$\log F'=L+WF/2-wz/2$ show that $\log F'$ has finite Dirichlet
integral on $V$: the additional term is controlled by the finite area
of $F(V)$. Put $A=\log F'$ and take smaller half-disks $V_r$ at
the source boundary point. Their quasidisk constants are uniform, and
$E_r=\int_{V_r}|A'|^2\to0$.

For $\delta(z)=\dist(z,\partial V_r)$, the analytic mean-value and
Cauchy estimates give
$\delta|A'|\le C\sqrt{E_r}$ and
$\delta^2|A''|\le C\sqrt{E_r}$.
Since $SF=A''-(A')^2/2$ and the hyperbolic density is bounded below
by a constant times $\delta^{-1}$, its hyperbolic Schwarzian norm
is at most $C(\sqrt{E_r}+E_r)$. Ahlfors' quasidisk extension criterion \cite{Ahlfors1963} (see also \cite[Theorem 4.1]{MVW}) gives local quasiconformal extensions for
small $r$; no prior image smoothness is required.

The $H^{1/2}$ trace of $A$ on a smaller boundary interval has
exponentials of its real part in every finite $L^p$. The Poisson
formula and Jensen's inequality give local Hardy $H^p$ bounds for
$F'=e^A$. The local Hardy boundary theorem therefore makes $F$
absolutely continuous there, with boundary derivative $e^A$ a.e.
This derivative is integrable and nonzero a.e. It follows that the
image arc is rectifiable and its arclength element is $|F'|\dd t$.
\end{proof}

\begin{remark}[Use of the two conclusions]
The arclength formula is used to compare the classical boundary
lengths on the two sides of a cut. It uses the finite-energy
estimates on $\log F'$ in the proof, not just the existence of a
quasiconformal extension. The extension provides coordinates for
Sobolev traces and gluing, as used in the conformal energy comparison
in Lemma~\ref{lem:weightedprojection}: if its distortion on $V$ is
at most $K$, then
\[
 \int_V|\nabla(h\circ\widetilde F)|^2\dd^2z
 \le K\int_{\widetilde F(V)}|\nabla h|^2\dd^2\zeta.
\]
In particular, matching traces can be glued after straightening the
image arc, and this arc has zero planar area. The separate conformal welding
argument in Theorem~\ref{thm:glue} constructs its charts from the
limiting boundary length correspondence.
\end{remark}

\begin{proof}[Proof of Theorem~\ref{thm:main}]
Let $\eta_\kappa$ be the strip interface. Under $q(z)=-e^{-z}$ its
half-plane image is $\mathrm{SLE}_\kappa(\rho_1;\rho_2)$, generated by
\[
 \partial_t g_t(z)=\frac{2}{g_t(z)-U_t},\qquad g_0(z)=z,
\]
with driver and force points satisfying \eqref{eq:SLE}.
The quantum welding law gives
\[
 \operatorname{Law}(A_\kappa,\eta_\kappa)
 =\operatorname{Law}(A_\kappa)\otimes
             \operatorname{Law}(\eta_\kappa).
\]
The field-dependent boundary-distance translation preserves this
independence: conditionally on the field it is a fixed strip
translation, which leaves the SLE law invariant.

Project Proposition~\ref{prop:joint} to the curve, obtaining a good
rate $J$. Recovering $G=(G_1,G_2)$ by \eqref{eq:derivedG} is
continuous at finite-rate points by Lemma~\ref{lem:trace}; contraction
therefore also gives a good rate $J_G$ for the geometric maps.
Since $G$ depends only on the curve, the product LDP identifies the
joint rate of the ambient data and geometric maps as
\[
 I(A)+J_G(G).
\]
Let $A_*$ be the unique zero of the ambient rate, given by
$\phi_W^*$ and its classical enhancements. If $J_G(G)<\infty$,
good-rate attainment in the contraction supplies a finite joint
preimage over $(A_*,G)$. Its interior cutting identities force
the cut fields \eqref{eq:cutzero}, up to the source translations
between geometric and field gauges. In this deterministic ambient
field both outer distance-one marks have strip coordinate zero.
Consequently \eqref{eq:derivedG} gives $s_1=s_2=0$ and $F_i=G_i$,
so the joint rate is exactly \eqref{eq:mainrate}.

Conversely let a Jordan crosscut have cut fields of finite energy.
Its outer maps are analytic on bounded intervals of the outer boundary by reflection;
ordinary change of variables gives their unit distance-one
constraints. Lemma~\ref{lem:localregularity} gives locally rectifiable
seam boundaries and local quasiconformal charts.
Taking Sobolev traces of the cutting identities in the local
quasiconformal charts yields the same classical length measure
$e^{\phi_W^*/2}\dd s$ from both sides. The cumulative lengths from
the finite end therefore give the same correspondence. These
profiles, their measures, and their cumulative functions form a
finite-rate component pair.

Proposition~\ref{prop:joint} supplies a finite limiting actual
configuration over this pair. The explicit deterministic cutting is
another configuration satisfying the hypotheses of
Theorem~\ref{thm:glue}. Its uniqueness, fixed by the outer mark, makes
them identical. This proves sufficiency and the rate formula for the
explicit deterministic cutting.

Lemma~\ref{lem:trace} makes the finite-rate trace and geometric maps
determine each other, and contraction infima are attained. Hence
$J(\eta)=J_G(G_1^\eta,G_2^\eta)$, proving \eqref{eq:mainrate} for
the Hausdorff trace as well. Both rates are infinite outside the
corresponding genuine curve and map configurations.
The LDPs restrict to the respective relative curve topologies: an
ambient closure of a relatively closed set adds no finite-rate
points, and relatively open sets are intersections with ambient
open sets. Compact sublevels remain compact in the curve space.

Finally the unique zero component pair is given by
$\phi_{W_i}^*$. Its explicit gluing into $\phi_W^*$ has
$F_1(z)=(W_1/W)z$ and
$F_2(z)=(W_2/W)z+i\pi W_1/W$. These maps satisfy all marks and field
identities. Hence the unique zero curve is the horizontal line
$\operatorname{Im}z=\pi W_1/W$, whose image under $q$ is the ray
of angle $\pi W_2/W$ in the upper half-plane.
\end{proof}

\section{The driving energy}
\label{sec:driving}

Write $q(z)=-e^{-z}$ and represent the curve in $\mathbb H$ by $q(\eta)$,
oriented from $0$ to $\infty$. Throughout this section $\eta$ denotes
this half-plane image and $J$ denotes the transported welding rate.
Its capacity parametrization is chosen so that
$\partial_tg_t(z)=2/(g_t(z)-U_t)$. Set
\[
 c=\frac{W-2}{2},\qquad
 V_t^-=g_t(0^-),\quad V_t^+=g_t(0^+).
\]
The last expressions denote the two prime-end images. For a Jordan
crosscut with no other real boundary contacts,
$V_t^-<U_t<V_t^+$ for $t>0$, and
$\dot V_t^\pm=2/(V_t^\pm-U_t)$ locally on $(0,\infty)$.
Put $L_t=U_t-V_t^-$, $R_t=V_t^+-U_t$, and $S_t=L_t+R_t$.
Define
\begin{equation}\label{eq:drivingenergy}
 \mathcal E_{\rho_1,\rho_2}(\eta)
 =\frac12\int_0^\infty
 \left|\dot U_t-\frac{\rho_1}{L_t}
                      +\frac{\rho_2}{R_t}\right|^2\dd t,
\end{equation}
with value $+\infty$ unless $U\in AC_{\rm loc}(0,\infty)$.
The integral is the monotone limit of its restrictions to
$[\varepsilon,T]$; it does not separately square the singular drift
and the driver at time zero.

\begin{theorem}[Equality of the two energies]\label{thm:driving}
For every fixed $W_1,W_2>0$ and every Jordan crosscut from $0$ to
$\infty$ in $\mathbb H$, with no other real boundary contacts,
\[
 J(\eta)=\mathcal E_{W_1-2,W_2-2}(\eta),
\]
including equality of the finite domains. In particular $J$ is ordinary
chordal Loewner energy when $W_1=W_2=2$.
\end{theorem}

For ordinary Loewner energy the conformal-map equivalence is due to
Wang~\cite{Wang}, whose proof uses additivity, an infinitesimal Loewner
calculation, and approximation. Here the logarithmic terms associated
with the force points change both the harmonic function and the
driving action. We prove the weighted comparison through Dirichlet
projections and the chordal Hadamard formula, as in the analytic
mechanism underlying SLE/free-field couplings~\cite{IK}. The comparison
keeps the endpoint cancellations together, rather than subtracting
separately divergent energies.

\subsection{The conformal comparison}

The comparison uses orthogonal projections in the zero-boundary
Dirichlet space. We fix their notation here; the isometry and the
boundary uniqueness argument are proved in
Appendices~\ref{app:projection} and~\ref{app:bloch}.

The Dirichlet space $\dot H_0^1(\mathbb H)$ is the completion of
$C_c^\infty(\mathbb H)$ with norm
$\|v\|_\nabla^2=\pi^{-1}\int_{\mathbb H}|\nabla v|^2$.
Spaces on subdomains are embedded by zero extension, or equivalently
by their quasi-continuous representatives vanishing quasi-everywhere
on the complement. For a zero-area simple Loewner slit $K_t$, set
\[
 O_t=\mathbb H\setminus K_t,\qquad
 \mathcal H_t=\dot H_0^1(\mathbb H)\ominus\dot H_0^1(O_t),\qquad
 P_t:\dot H_0^1(\mathbb H)\longrightarrow\mathcal H_t.
\]
Thus $\mathcal H_t$ consists of the functions harmonic on $O_t$.
The projections are increasing in $t$. Define
\[
 p_t(z)=-\Ima\frac1{g_t(z)-U_t}>0,\qquad z\in O_t.
\]

Let $O_i$ be the two components of $\mathbb H\setminus\eta$,
with $O_1$ adjacent to $\R_-$, and put
\[
 h_i=q\circ G_i^{-1}\circ q^{-1}:O_i\longrightarrow\mathbb H.
\]
Positive dilations of $h_i$ do not affect the expressions below.
Choose holomorphic logarithms and define
\begin{equation}\label{eq:weightedA}
 A_i(z)=\log h_i'(z)+\frac{\rho_i}{2}\log h_i(z)-c\log z.
\end{equation}
Add imaginary constants to make $v_i=\Ima A_i$ zero on its outer
real arc. The change of coordinates can be checked explicitly:
\[
 \Rea A_i(q(G_i(\zeta)))
 =-\log|G_i'(\zeta)|+\frac W2\Rea G_i(\zeta)
                     -\frac{W_i}{2}\Rea\zeta+\text{constant}.
\]
Thus $\phi_i^\eta+W_i\Rea\zeta$ is $-2\Rea A_i\circ q\circ G_i$
up to an additive constant. Conformal invariance and equality of
the Dirichlet energies of the real and imaginary parts give
\begin{equation}\label{eq:Anorm}
 J(\eta)=\frac1\pi\sum_i\int_{O_i}|\nabla v_i|^2\dd^2z.
\end{equation}

\begin{lemma}[Matching and harmonic projection]\label{lem:weightedprojection}
Suppose $J(\eta)<\infty$. The functions $v_i$ glue to
$v\in\dot H_0^1(\mathbb H)$, harmonic off $\eta$, with
$\|v\|_\nabla^2=J(\eta)$. Define on $O_t$
\begin{align}
 H_t(z)={}&\log g_t'(z)+\log(g_t(z)-U_t)
       +\frac{\rho_1}{2}\log(g_t(z)-V_t^-)
       \notag\\[-2pt]
 &+\frac{\rho_2}{2}\log(g_t(z)-V_t^+)-c\log z.
 \label{eq:weightedHt}
\end{align}
The branches are continued from infinity. Then
\begin{equation}\label{eq:Htprojection}
 \Ima H_t=P_tv.
\end{equation}
No absolute continuity of $U$ is assumed in this assertion.
\end{lemma}
\begin{proof}
Lemma~\ref{lem:localregularity} gives local quasiconformal charts,
rectifiability, and logarithmic-derivative $H^{1/2}$ traces on finite
seam arcs. Let $\beta$ denote the tangent-angle trace oriented from
$0$ to $\infty$. On the seam $h_1>0$ and $h_2<0$. The boundary
derivatives have arguments $-\beta$ and $\pi-\beta$, respectively.
On the outer arcs, $h_1<0$ and $h_2>0$, with positive derivatives.
Consequently the two imaginary traces, after the stated constant
choices, both equal
\begin{equation}\label{eq:weightedwinding}
 \frac{\pi W_2}{2}-\beta(z)-c\arg z.
\end{equation}
The logarithm branches are fixed from the outer arcs. The tangent
angle is unwrapped consistently; Lemma~\ref{lem:ends} fixes the
relative branch at an end, so the displayed agreement has no
undetermined multiple of $2\pi$. Equivalently, the difference is a
locally $H^{1/2}$ integer-valued function and is constant on the
connected seam, with its constant fixed by those end branches.

Sobolev gluing in the local quasiconformal charts gives a locally
$H^1$ function across the seam. Its total Dirichlet energy is
\eqref{eq:Anorm} and it has zero trace on $\R\setminus\{0\}$.
The exceptional endpoints have zero Sobolev capacity. Zero extension
across $\R$, or the half-plane Hardy inequality and exhaustion,
therefore puts this function in $\dot H_0^1(\mathbb H)$. This proves the first
assertion. In particular the slit has zero area, by its local
quasiconformality away from its endpoints.

On the real boundary $\Ima H_t=0$. On the two faces of the past
slit, the arguments of $g_t-U_t$ are $\pi$ and $0$; the arguments
of $g_t-V_t^-$ and $g_t-V_t^+$ are $0$ and $\pi$ on both faces.
Together with the boundary argument of $g_t'$, this gives exactly
\eqref{eq:weightedwinding}. Thus $\Ima H_t$ has the same boundary
trace as $P_tv$ away from $0$, the tip, and infinity.

Here the trace assertion is distributional, not just pointwise a.e.
Indeed, near a past seam point different from the tip and the
origin, $g_t\circ h_i^{-1}$ maps a real interval homeomorphically
to a real interval. Schwarz reflection extends it analytically
with nonzero derivative. On that neighborhood $H_t-A_i$ is an
analytic function of $h_i$ with bounded derivative on a smaller
half-disk. Its local Dirichlet integral is finite by conformal
invariance. Thus both functions have local $H^1$ traces there.
The outer real intervals are simpler reflection intervals.

Put $F_t=g_t^{-1}$. The pullback $H_t\circ F_t$ is
\[
 -\log F_t'+\log(w-U_t)
 +\frac{\rho_1}{2}\log(w-V_t^-)
 +\frac{\rho_2}{2}\log(w-V_t^+)-c\log F_t.
\]
It is Bloch, without any energy assertion at the marked points.
The pullback of $P_tv$ is Bloch because it has finite Dirichlet
energy. Their imaginary-part difference has zero distributional
trace outside $V_t^-,U_t,V_t^+$ and infinity.
Lemma~\ref{lem:blochunique} proves \eqref{eq:Htprojection}.
In particular it establishes the finite energy of $H_t$, rather
than presupposing it.
\end{proof}

\begin{proof}[Proof of Theorem~\ref{thm:driving}]
\emph{From conformal energy to driving energy.}
Let $J(\eta)<\infty$. Lemmas~\ref{lem:projection} and
\ref{lem:weightedprojection} give a coefficient $e$ such that
\[
 \Ima H_t(z)=\int_0^t p_s(z)e(s)\dd s.
\]
Fix $0<\varepsilon<T$ and $z\in O_T$. The right side is absolutely
continuous in $t$. In \eqref{eq:weightedHt} all terms other than
$\arg(g_t(z)-U_t)$ are already absolutely continuous on
$[\varepsilon,T]$: this uses the ordinary Loewner equations and
the positive force gaps, not regularity of $U$. Hence
$\theta_t=\arg(g_t(z)-U_t)$ is absolutely continuous. Write
$g_t(z)=x_t+iy_t$. On the compact time interval $\theta_t$ stays
in a compact subset of $(0,\pi)$, and
\[
 U_t=x_t-y_t\cot\theta_t.
\]
Since $x_t,y_t$ are absolutely continuous, this proves
$U\in AC[\varepsilon,T]$.

We may now differentiate. The Loewner equations give
\[
 \begin{aligned}
 \partial_t\log g_t'(z)&=-\frac{2}{(g_t(z)-U_t)^2},\\
 \partial_t\log(g_t(z)-U_t)
   &=\frac{2}{(g_t(z)-U_t)^2}
                  -\frac{\dot U_t}{g_t(z)-U_t},\\
 \partial_t\log(g_t(z)-V_t^\pm)
   &=\frac{2}{(g_t(z)-U_t)(U_t-V_t^\pm)}.
 \end{aligned}
\]
The squared-denominator terms cancel. Substituting
$U_t-V_t^-=L_t$ and $U_t-V_t^+=-R_t$ yields
\begin{equation}\label{eq:Htderivative}
 \partial_tH_t(z)=
 -\frac{\dot U_t-\rho_1/L_t+\rho_2/R_t}{g_t(z)-U_t}.
\end{equation}
Comparison of imaginary parts, using $p_t(z)>0$, identifies
$e(t)=\dot U_t-\rho_1/L_t+\rho_2/R_t$ a.e.
Finally $P_tv\to v$ in Dirichlet norm, so
\eqref{eq:projectionisometry} gives $J(\eta)=\mathcal E_{\rho_1,\rho_2}(\eta)$.
This implication has derived both absolute continuity and finiteness
of the driving action from the conformal energy.
\medskip
\noindent\emph{From driving energy to conformal energy.}
Suppose instead that \eqref{eq:drivingenergy} is finite, and write
$e(t)=\dot U_t-\rho_1/L_t+\rho_2/R_t$. On every
$[\varepsilon,T]$ the drift is bounded, so $\dot U\in L^2$ there.
Finite ordinary Loewner-energy regularity \cite{WangLDP,Wang},
applied after mapping out an earlier initial segment, shows that every compact
interior subarc has zero area. A countable union shows that the
whole slit has zero area; this use of ordinary finite energy is
only in this converse implication.

Lemma~\ref{lem:projection} constructs compatible vectors $v_t$ from
$e|_{[0,t]}$. Their increments are orthogonal, so they converge in
$\dot H_0^1(\mathbb H)$ to $v$ with
$\|v\|_\nabla^2=\tfrac12\int_0^\infty e^2$ and $P_tv=v_t$.
This limit is harmonic off the full slit, since all the $v_t$ are
harmonic there and Dirichlet convergence preserves distributional
harmonicity.
Equation~\eqref{eq:Htderivative} holds on compact positive time
intervals. At each fixed $z\in\mathbb H$,
$H_t(z)\to0$ as $t\downarrow0$, since $U_t,V_t^\pm\to0$,
$g_t\to\operatorname{id}$, and $1+(\rho_1+\rho_2)/2=c$.
Also $p_s(z)$ is bounded near zero and $e\in L^2$.
Integrating therefore gives $\Ima H_t=v_t$ on $O_t$.

Define $A_i$ by \eqref{eq:weightedA}, without assuming its energy
finite. Every compact subarc of the seam away from its endpoints lies in the past of some time $t$.
On that arc Schwarz reflection of $g_t\circ h_i^{-1}$ again shows
that $H_t-A_i$ is locally an analytic function of $h_i$.
Since $H_t$ now has finite energy, this proves local $H^1$ traces
for $A_i$ and their agreement with $v$ after the same imaginary
constant choices. On the outer arcs the traces are zero.

Writing $f_i=h_i^{-1}$, the function
$A_i\circ f_i=-\log f_i'+(\rho_i/2)\log w-c\log f_i$
is Bloch. The pullback of $v|_{O_i}$ has finite energy and is
harmonic. Their imaginary-part difference has zero distributional
trace outside the two end points. Lemma~\ref{lem:blochunique}
identifies them. Thus both $A_i$ have finite Dirichlet energy,
and \eqref{eq:Anorm} proves $J(\eta)=\|v\|_\nabla^2$.
The boundary-distance gauge constraints hold by change of variables
on the analytic outer arcs, as in the converse part of
Theorem~\ref{thm:main}. This proves equality of the finite domains.
\end{proof}

\subsection{Finite-cutoff boundary terms and the endpoints}

The projection cutoff gives a precise version of the weighted
boundary calculation. Let $g_{\varepsilon,T}=g_T\circ g_\varepsilon^{-1}$,
and put $D_\pm=g_{\varepsilon,T}'(V_\varepsilon^\pm)>0$. Define
\begin{align}
 B_{\varepsilon,T}={}&
 -\rho_1\log\frac{L_T}{L_\varepsilon}
 -\rho_2\log\frac{R_T}{R_\varepsilon}
 -\frac{\rho_1(\rho_1+4)}4\log D_-
 \notag\\
 &-\frac{\rho_2(\rho_2+4)}4\log D_+
 -\frac{\rho_1\rho_2}{2}\log\frac{S_T}{S_\varepsilon}.
 \label{eq:boundarycorrection}
\end{align}
\begin{corollary}[Combined boundary cancellation]\label{cor:boundarycancel}
For a finite-rate curve,
\begin{equation}\label{eq:finitecomparison}
 \|(P_T-P_\varepsilon)v\|_\nabla^2
 =\frac12\int_\varepsilon^T\dot U_t^2\dd t+B_{\varepsilon,T}.
\end{equation}
The combined right side tends to $J(\eta)$ as
$\varepsilon\downarrow0$, $T\uparrow\infty$.
\end{corollary}
\begin{proof}
We have already proved local absolute continuity and finite driving
action without using this expansion. On a finite interval,
\begin{gather*}
 \dot L=\dot U+2/L,\qquad \dot R=2/R-\dot U,\qquad
 (\log S)'=2/(LR),\\
 \log D_-=-2\int_\varepsilon^TL^{-2}\dd t,\qquad
 \log D_+=-2\int_\varepsilon^TR^{-2}\dd t.
\end{gather*}
Expanding $\tfrac12(\dot U-\rho_1/L+\rho_2/R)^2$ and
integrating these exact logarithmic derivatives gives
\eqref{eq:boundarycorrection}. Orthogonality in
Lemma~\ref{lem:projection} proves \eqref{eq:finitecomparison}.
Its limit is $\|v\|_\nabla^2$, since $P_Tv\to v$ and
$\|P_\varepsilon v\|_\nabla^2=\tfrac12\int_0^\varepsilon e^2\to0$.
No individual divergent term is assigned a finite limit.
\end{proof}

For completeness there is an explicit positive expression for the
same cancellation. Set $s=\log S_t$, $r(s)=L_t/S_t$,
$p=W_1/W$, and
\[
 b(r)=(W_1-Wr)/2,\qquad
 V(r)=2W_1\log(p/r)+2W_2\log((1-p)/(1-r)).
\]
Then $\dd t/\dd s=e^{2s}r(1-r)/2$, and direct substitution gives,
for $\alpha=\log S_\varepsilon$, $\beta=\log S_T$,
\begin{align}
 \frac12\int_\varepsilon^T\dot U_t^2\dd t+B_{\varepsilon,T}
 &=\int_\alpha^\beta\frac{|\dot r-b(r)|^2}{r(1-r)}\dd s
 \notag\\
 &=\int_\alpha^\beta
 \frac{\dot r^2+(W^2/4)(r-p)^2}{r(1-r)}\dd s
 +\frac{V(r(\beta))-V(r(\alpha))}{2}.
 \label{eq:positivecancellation}
\end{align}
Since the whole driving action is now known finite,
$u=\dot r-b(r)\in L^2(\R)$. Boundedness of $r$ and the linear
equation $\dot r+(W/2)(r-p)=u$ give
\[
 r(s)-p=\int_{-\infty}^s e^{-W(s-v)/2}u(v)\dd v.
\]
Thus $r\to p$ at both ends and $V(r)\to0$. This removes the
last term in \eqref{eq:positivecancellation} rigorously. Both
$W_i>0$ are used here; they put $p$ strictly inside $(0,1)$.
Finite action need not make $r-p$ integrable, so an individual
$\log D_\pm$ can retain slower divergent corrections after its
leading linear term is removed. The combined identities above
are the appropriate endpoint renormalization.

\clearpage
\appendix

\section{Dirichlet projections along a Loewner chain}\label{app:projection}

We retain the notation of Section~\ref{sec:driving}. The following
isometry is deterministic. Hadamard-based isometries between $L^2$ spaces and Dirichlet
spaces appear in \cite{HedenmalmNieminen2014} and, in the
Loewner--Kufarev setting, in \cite{ViklundWang2024}.
We give the chordal projection construction needed here
for continuous driving functions. In particular it requires no driving
LDP and no absolute continuity of the Loewner driver.

\begin{lemma}[Projection isometry]\label{lem:projection}
Assume that the simple Loewner slit $K_T$ has zero planar area.
For $T<\infty$ there is a unitary map
$C_T:L^2(0,T)\longrightarrow\mathcal H_T$ such that, for $z\in O_T$,
\[
 (C_T f)(z)=\sqrt2\int_0^T p_t(z)f(t)\dd t,
 \qquad P_s C_T f=C_T(\mathbf1_{[0,s]}f),\quad 0\le s\le T.
\]
Consequently, for $v\in\dot H_0^1(\mathbb H)$, there is a unique compatible
coefficient $e\in L^2_{\rm loc}(0,\infty)$ satisfying
\begin{equation}\label{eq:projectionisometry}
 P_tv(z)=\int_0^t p_s(z)e(s)\dd s,\qquad
 \|P_tv\|_\nabla^2=\frac12\int_0^t|e(s)|^2\dd s.
\end{equation}
If $v$ is harmonic off $\bigcup_tK_t$, then $P_tv\to v$ in
$\dot H_0^1(\mathbb H)$ as $t\to\infty$.
\end{lemma}
\begin{proof}
Normalize the Dirichlet Green function by $-\Delta G=2\pi\delta$.
For $z,w\in O_T$, Green's identity shows that the reproducing kernel
of $\mathcal H_t$ is
\[
 k_t(z,w)=\tfrac12\{G_{\mathbb H}(z,w)-G_{O_t}(z,w)\}.
\]
The factor $1/2$ comes from our norm $\pi^{-1}\int|\nabla v|^2$.
The logarithmic poles cancel, so these kernel functions have finite
Dirichlet energy. Their linear span is dense in $\mathcal H_T$:
a vector orthogonal to all of them vanishes on $O_T$, and hence
vanishes a.e., since $K_T$ has zero planar area.

The formula
\[
 G_{O_t}(z,w)=\log\left|
 \frac{g_t(z)-\overline{g_t(w)}}{g_t(z)-g_t(w)}\right|
\]
and the Loewner equation, requiring only continuity of $U$, give
\begin{equation}\label{eq:hadamard}
 -\partial_tG_{O_t}(z,w)=4p_t(z)p_t(w),\qquad
 k_t(z,w)=2\int_0^t p_s(z)p_s(w)\dd s.
\end{equation}
Kernel factorization initially gives a coisometry $C_T$ onto
$\mathcal H_T$, with the stated evaluation formula. It remains to
prove injectivity; kernel factorization alone would give only a
minimum-norm representation.

Let $M_s$ be multiplication by $\mathbf1_{[0,s]}$ on $L^2(0,T)$.
Equation~\eqref{eq:hadamard} implies
$C_T M_s C_T^*=P_s|_{\mathcal H_T}$. Since $C_TC_T^*=1$,
$C_T^*$ is an isometry. For $h\in\mathcal H_T$,
\[
 \|M_sC_T^*h\|_2^2=\langle P_sh,h\rangle_\nabla
 =\|P_sh\|_\nabla^2=\|C_TM_sC_T^*h\|_\nabla^2.
\]
Equality of these norms puts $M_sC_T^*h$ in
$(\ker C_T)^\perp=\operatorname{ran}C_T^*$, and gives
$M_sC_T^*h=C_T^*P_sh$.
Thus the closed range of $C_T^*$ is invariant under all $M_s$,
and, by self-adjointness, reduces them. A closed reducing subspace
for all interval multipliers is $L^2(E)$ for a measurable
$E\subset(0,T)$. But it contains
$C_T^*k_T(z,\cdot)=\sqrt2p_\cdot(z)$, which is strictly positive
on $(0,T)$ for any fixed $z\in O_T$. Hence $E=(0,T)$ a.e.
This proves unitarity and the intertwining formula.

Apply $C_T^{-1}$ to $P_Tv$ and set $e=\sqrt2C_T^{-1}P_Tv$.
Intertwining makes these coefficients consistent as $T$ varies,
and proves \eqref{eq:projectionisometry}. Finally, the projections
increase to the orthogonal complement of
$\bigcap_t\dot H_0^1(O_t)$. The quasi-continuous characterization
identifies this intersection with the functions vanishing on the
whole slit. A function harmonic off that slit belongs to the
orthogonal complement, proving the last assertion.
\end{proof}

\clearpage
\section{Boundary uniqueness}\label{app:bloch}

This appendix supplies the uniqueness statement used to identify
the charged harmonic function with its Dirichlet projection.
Its boundary hypothesis is distributional, so it also excludes
singular contributions supported at marked prime ends.

\begin{lemma}[Bloch uniqueness at finitely many points]
\label{lem:blochunique}
Let $u$ be harmonic in the disk with
$\sup_{z\in\D}(1-|z|^2)|\nabla u(z)|<\infty$.
If its distributional boundary trace vanishes on the complement
of a finite subset of $\partial\D$, then $u=0$.
The corresponding assertion holds in $\mathbb H$, with infinity
included among the possible exceptional points.
\end{lemma}
\begin{proof}
Integrating the gradient radially gives
$|u(re^{i\theta})|\le C+C\log(1/(1-r))$, uniformly in $\theta$.
For example, estimating the $n$th Fourier coefficient on the circle
$r=1-1/n$ bounds its boundary coefficient by $O(\log(n+1))$.
Thus $u$ has a boundary distribution and is its Poisson extension.
A distribution supported at finitely many points is a
finite sum of derivatives of point masses. At a point supporting a
nonzero term, its Poisson extension grows like a positive power of
$(1-r)^{-1}$, as seen along radial or fixed non-tangential approaches.
Terms at the other points remain bounded there; the highest nonzero
derivative cannot cancel. This contradicts the logarithmic growth.
Hence the distribution, and then $u$, is zero. A M\"obius map gives
the half-plane statement.
\end{proof}

We use the distributional condition, not just a.e. boundary values.
In the applications it follows from local $H^1$ traces on every
boundary interval away from the explicitly listed marked points.
For an injective holomorphic map $F:\mathbb H\to\mathbb H$, both
$\log F'$ and $\log F$ are analytic Bloch functions: use the
pre-Schwarzian estimate for the former and the bounded imaginary
part of $\log F$ for the latter. A harmonic function of finite
Dirichlet energy is Bloch by the analytic mean-value estimate.
These facts allow Lemma~\ref{lem:blochunique} to establish finite
energy of a candidate by comparison with a projection; finite
energy of that candidate is not a hypothesis.

\end{document}